\documentclass[nonblindrev]{informs4}

\OneAndAHalfSpacedXI

\usepackage{amsmath,amssymb,amsfonts}

\usepackage{natbib}
 \bibpunct[, ]{(}{)}{,}{a}{}{,}%
 \def\bibfont{\small}%

\usepackage{rotating}
\usepackage{fancyvrb}
\usepackage{appendix}
\usepackage{graphicx}
\usepackage{hyperref}
\usepackage{subcaption}
\usepackage{url}

\TheoremsNumberedThrough     
\ECRepeatTheorems
\JOURNAL{ }

\EquationsNumberedThrough    

\newcommand{\E}{\mathbb{E}}
\newcommand{\R}{\mathbb{R}}
\newcommand{\Prob}{\mathbb{P}}
\newcommand{\X}{\mathcal{X}}
\newcommand{\Px}{\mathcal{P}}

\newtheorem{coro}{Corollary}
\gdef\AQ#1{}
\gdef\CQ#1{}

\begin{document}


	\RUNAUTHOR{Moucer~et~al.} %

	\RUNTITLE{Constructive approaches to concentration inequalities}

\TITLE{Constructive approaches to concentration inequalities with independent random variables}

	\ARTICLEAUTHORS{

\AUTHOR{Céline Moucer}
\AFF{Inria, Département d'Informatique de l'Ecole Normale Supérieure, PSL Research University, Paris, France. \url{celine.moucer@inria.fr}}
\AFF{Ecole Nationale des Ponts et Chaussées, Marne-la-Vallée, France.}

\AUTHOR{Adrien Taylor}
\AFF{Inria, Département d'Informatique de l'Ecole Normale Supérieure, PSL Research University, Paris, France. \url{adrien.taylor@inria.fr}}

\AUTHOR{Francis Bach}
\AFF{Inria, Département d'Informatique de l'Ecole Normale Supérieure, PSL Research University, Paris, France. \url{francis.bach@inria.fr}}
}

	

	\ABSTRACT{Concentration inequalities, a major tool in probability theory, quantify how much a random variable deviates from a certain quantity. This paper proposes a systematic convex optimization approach to studying and generating concentration inequalities with independent random variables. Specifically, we extend the generalized problem of moments to independent random variables. 
 
 We first introduce a variational approach that extends classical moment-generating functions, focusing particularly on first-order moment conditions. Second, we develop a polynomial approach, based on a hierarchy of sum-of-square approximations, to extend these techniques to higher-moment conditions. Building on these advancements, we refine Hoeffding's, Bennett's and Bernstein's inequalities, providing improved worst-case guarantees compared to existing results.}
 

\SUBJECTCLASS{60E15, 90C22, 90C25, 60F10}


	
\maketitle
	
\section*{Introduction.}

Concentration inequalities have emerged as a major tool in probability theory, finding applications in learning theory~\citep{1996_devroye}, in random matrix theory~\citep{2012Tao} or statistical physics or mechanics~\citep{1998_Dembo}. These inequalities quantify how much a random variable deviates from a certain quantity, usually its mean. Classical examples include Markov's inequality for bounding probabilities of deviations from zero or Chebyshev's inequalities for deviation from the mean. In machine learning and statistics, where the data are often assumed to be independent and identically distributed (i.i.d.), basic inequalities such as Hoeffding's inequality, Bennett's inequality or Bernstein's inequality~\cite[Chapter 6]{2013Boucheron} are extensively used, e.g., for characterizing generalization properties of machine learning algorithms~\citep{2024Bach}. For instance, Hoeffding's inequality~\citep{1963Hoeffding} states that for $X_1, \ldots, X_n$ i.i.d. random variables taking their values almost surely in $[0, 1]$, the sum $\sum_{i=1}^n X_i$ has a subgaussian tail with deviation $t \geqslant 0$:
\begin{equation*}
    \Prob\left(\sum_{i=1}^n (X_i - \E[X_i]) \geqslant  nt)\right) \leqslant \exp\left( -2 n t^2\right).
\end{equation*}

In this work, we provide a principled approach to concentration inequalities for independent univariate random variables with finite moments. Let $\X$ be a subset of $\R$ and $\Px(\X)$ be the set of distributions on~$\X$. Given $X_1, \ldots, X_n$ independent random variables generated from distributions $p_1, \ldots, p_n \in \Px(\X)$, we formulate the generalized problem of moments for independent random variables as follows:
\begin{equation}
\label{eq:independent_upper_bound_problem}
    \begin{aligned}
    \rho_n &= \sup_{p_1, \ldots, p_n \in \Px(\X)} \E_{p_1,\ldots, p_n}[F(X_1, \ldots, X_n)] \ {\rm \ such \ that \ } \forall i, \E_{p_i}[g_i(X_i)] = \mu_i, \\
 &= \sup_{p_1, \ldots, p_n \in \Px(\X)} \int_{\X^n} F(x_1, \ldots, x_n) dp_1(x_1)\cdots dp_n(x_n) \ {\rm \ such \ that \ } \forall i, \int_{\X} g_i(x_i) dp_i(x_i) = \mu_i,
\end{aligned}
\end{equation}
for some functions $F:x \in \X^n \mapsto \R^+$ and $g_i:x \in \X \mapsto \R^m$. For instance, Hoeffding's inequality involves $F(x) = \mathbf{1}_{\sum_{i=1}^n x_i \geqslant nt + \sum_{i=1}^n \E[X_i]}$, and therefore $\E_{p_1, \ldots, p_n}[F(x)] = \Prob(\sum_{i=1}^n (X_i - \E[X_i]) \geqslant nt)$. Problem~\eqref{eq:independent_upper_bound_problem} is an infinite-dimensional non-convex problem. Without further assumptions on $F$ and $g$, minimizing with respect to distributions $p_i$ is often intractable.

This problem is closely related to the generalized problem of moments formalized by \cite{2008_Lasserre}, which extends the traditional problem of moments~\citep{1987_Landau} that seeks a measure matching a given set of moments. The search for optimal multivariate Chebyshev's inequalities began in the 1960s~\citep{Marshall1960, 1962_Isii}. \cite{1962_Isii, 1964_Isii_InequalitiesOT}, along with \cite{1966Karlin_Studden}, formalized the pursuit of sharp inequalities by framing it as an optimization problem. Let $\mathcal{Y} \subset \R^n$ and $\Px(\mathcal{Y})$ be the set of probability distributions on $\mathcal{Y}$. The generalized problem of moments takes the form:
\begin{equation}
\label{eq:multivariate_upper_bound}
\begin{aligned}
  \rho &=   \sup_{p \in \Px(\mathcal{Y} )} \ \E_p[F(X)], \ 
  {\rm such \ that \ } \E_p[g(X)]= \mu \\
    &= \sup_{p \in \Px(\mathcal{Y} )} \int_{\mathcal{Y} } F(x) dp(x) \
    {\rm \ such \ that \ } \int_{\mathcal{Y} }g(x)dp(x) = \mu,
\end{aligned}
\end{equation}
where $g:\mathcal{Y}  \mapsto \R^m$. Compared to Problem~\eqref{eq:independent_upper_bound_problem}, the generalized problem of moments optimizes over distributions $p \in \Px(\mathcal{Y} )$ that are not necessarily products of their marginals. Under mild assumptions on the moment vector~$\mu$, strong duality holds~\cite[Theorem 3.1]{1964_Isii_InequalitiesOT}. Its Lagrangian relaxation was first formulated by~\cite{1964_Isii_InequalitiesOT} as follows:
\begin{equation}
\label{eq:dual_generalized_p_moment}
\begin{aligned}
    \rho &= \inf_{\alpha \in \R, \beta \in \R^m} \alpha + \beta \mu \
    {\rm \ such \ that \ } \forall x \in \mathcal{Y} , F(x) \leqslant \alpha + \beta^\top g(x), \\
    &= \inf_{\alpha \in \R, \beta \in \R^m} \alpha + \beta \mu + \sup_{x \in \mathcal{Y}} \{ F(x) - (\alpha + \beta^\top g(x)) \}.
\end{aligned}
\end{equation}
where $\alpha \in \R$ corresponds to the dual variable associated to the constraint $\int_{\mathcal{Y} }dp(x) = 1$, and $\beta \in \R^m$ to the constraint $\int_{\mathcal{Y} } g(x) dp(x) = \mu$. The dual Problem~\eqref{eq:dual_generalized_p_moment} is convex as it is expressed as the pointwise supremum of affine functions. Yet, it is unclear how to deal with the constraint ``$\forall x \in \mathcal{Y} , F(x) \leqslant \alpha + \beta^\top g(x)$'', because it corresponds to infinitly many linear constraints in $(\alpha, \beta)$. 

Those problems are now traditionally approached via convex reformulations or approximations using semidefinite programming (SDP). \cite{2005Bertsimas_popescu} first investigated optimal bounds for $\E_p[F(X)] = \Prob(X \in S)$ assuming $\X$ and $S$ to be semi-algebraic sets. For univariate random variables, they efficiently solved it using a single SDP, allowing them to derive tight bounds. For multivariate random variables, they proposed a series of semidefinite relaxations using sum-of-square representations (SoS). More generally, \cite{2008_Lasserre} investigated the generalized problem of moments and derived a hierarchy of SDPs converging to the optimal value (extending the methodology developed for approximating global optimization problems~\citep{2001Lasserre}). Simultaneously, \cite{2007_Vandenberghe_Boyd_Comanor, Comanor2006} reformulated the generalized Chebyshev inequality as linear matrix inequalities (LMI) using an $S$-procedure. 

These SDP-based approaches, along with tight guarantees, often allow reconstructing corresponding worst-case distributions~\cite[Section 5.1]{2005Bertsimas_popescu}~\cite[Section 2.2]{2007_Vandenberghe_Boyd_Comanor}. These extremal distributions turns out to be discrete~\cite[Theorem 1]{Rogosinski1958} and may even be specified with $m+2$ Dirac, where $m$ is the number of constraints in~\eqref{eq:multivariate_upper_bound}. When the distributions in the problem under consideration are continuous with additional properties like symmetry or unimodality, these bounds may no longer be sharp. Therefore, \cite{2005_Popescu} generalized Chebyshev's inequality to convex classes of distributions generated by an appropriate parametric family of distributions. Building on Choquet's theory and conic duality, they provided a SDP reformulation (resp.~approximation) of the generalized problem of moments for such univariate (resp.~multivariate) distributions. From this framework, \cite{2015_VanParys_Goulart_Kuhn} extended Gauss inequalities to multivariate unimodal distributions and outlined a methodology for computing worst-case unimodal distributions related to this problem.

\paragraph{Research questions and assumptions.} Throughout this work, we assume $X_1, \ldots, X_n$ to be independent random variables in $\X \subset \R$ with finite moments $\forall i, \E_{p_i}[g(X_i)] = \mu_i \in \R^m$. For all $i$, we define $\mathcal{P}_{\mu_i}(\X) = \{p_i \in \mathcal{P}(\X), \int_{\X}g_i(x_i)dp_i(x_i) = \mu_i\}$ the set of (univariate) distributions on $\X$ with moment $\mu_i$. If in addition, $X_1, \ldots, X_n$ follows the same distributions with moments $\mu_1 = \mu_2 = \cdots = \mu_n$, they are said to be independent and identically distributed (i.i.d.). We assume that $F$ is an indicator function, that is for $S$ an appropriately selected compact semi-algebraic subset of $\X^n$ that emerges from the problem under consideration, $\forall x \in \X^n, F(x) = \mathbf{1}_{x \in S}$. The generalized moment problem for independent variables takes the form:
\begin{equation}
\label{eq:research_question}
    \begin{aligned}
        \rho_n &= \sup_{\forall i, p_i \in \Px(\X)} \int_{x \in \X^n} \mathbf{1}_{x \in S} dp_1(x_1) \cdots dp_n(x_n) \ {\rm  \ such \ that \ } \forall i, \ \int_{\X} g(x_i) dp_i(x_i) = \mu_i, \\
        &= \sup_{\forall i, p_i \in \Px_{\mu_i}(\X)} \int_{x \in \X^n} \mathbf{1}_{x \in S} dp_1(x_1) \cdots dp_n(x_n).
    \end{aligned}
\end{equation}
In this formulation, independence is encoded by explicitly accounting for the constraint ``$\forall x \in \X^n, p(x) = p_1(x_1)\cdots p_n(x_n)$''. Classical concentration inequalities correspond to upper bounds to these problems for specific choices of $g_i$'s. It is natural to wonder how tight such bounds are. To this end, we propose constructive and tractable approaches to upper bounding Problem~\eqref{eq:research_question}, along with a comparison to existing bounds. Then, we raise the issue of reconstructing worst-case distributions that can potentially match these bounds in some scenarios.

\paragraph{Contributions.} We propose constructive approaches for computing concentration inequalities of independent variables given a set of moments. To this end, we start by considering variational formulations to Problem~\eqref{eq:research_question}: 
\begin{equation}
\label{eq:research_question_upper_bound}
   \rho_n^{\mathcal{H}} = \inf_{H \in \mathcal{H}} \sup_{\forall i, p_i \in \mathcal{P}(\mu_i)} \int_{\X^n} H(x) dp_1(x) \cdots dp_n(x) {\rm \ such \ that \ } \forall x \in \X^n, F(x) \leqslant H(x),
\end{equation}
where $\mathcal{H}$ is a well-chosen set of functions. It is straightforward that Problem~\eqref{eq:research_question_upper_bound} yields an upper bound to the generalized problem of moments for independent random variables, that is:
\begin{equation*}
    \rho_n \leqslant \rho_n^{\mathcal{H}}.
\end{equation*} 

Efficiently leveraging the variational formulation~\eqref{eq:research_question_upper_bound} requires strategies to enforce the constraint $\forall x \in \X^n, F(x) \leqslant H(x)$ and to bound the objective. The choice of the function families $\mathcal{H}$ ensure these two requirements are satisfied. We propose two natural and complementary strategies based on convex optimization, each relying on different choices for $\mathcal{H}$.

\begin{enumerate}
    \item The first strategy revolves around a family of product-functions $\forall x \in \X^n, U(x) = \prod_{i=1}^n u_i(x_i)$ inspired from classical probability proofs and variational probabilistic inference~\citep{1999JaakkolaJordan}. Given such a function, we define a \textit{separable approach} by formulating Problem~\eqref{eq:research_question_upper_bound} as $n$ univariate subproblems. We then develop a \textit{variational approach} by optimizing over the family of product-funtions. When $F$ is log-convex and first-order moments are finite with $\forall i, \E[g_i(X_i)] = \E[X_i] = \mu_i$, the variational approach takes the form of a finite-dimensional convex optimization reformulation, that can be efficiently solved. This strategy significantly improves Hoeffding's inequality for small values of $n$ and accurately meets the asymptotic large deviations for large $n$. Moreover, this approach allows the reconstruction of distributions involved at the optimum of the upper bound Problem~\eqref{eq:research_question_upper_bound}. However, when it comes to higher-order finite moments, the variational approach formulates as a nonconvex problem that cannot be solved efficiently anymore.
    \item The second strategy relies on a family of polynomial upper bounds, which is particularly suited to higher-moment conditions. Referred to as a \textit{polynomial approach}, this strategy formulates as a non-convex optimization problem that can be effectively approximated by a series of sum-of-square formulations~\citep{2008_Lasserre}. This approach contributes to refining Bernstein's and Bennett's inequality, which are fundamental probabilistic bounds.
    \item Finally, we extend the polynomial approach to a \textit{feature-based approach} relying on broader families of upper bounds~$\mathcal{H}$. Compared to the variational approach, this method refines Hoeffding's inequality using higher-order polynomials when applied to two random variables. While this methodology introduces finer approximations $\rho_n \leqslant \rho_n^{\mathcal{H}}$, it often requires well-chosen relaxations to approximate~$\rho_n^{\mathcal{H}}$.
\end{enumerate}

\paragraph{Outline of the paper. } Section~\ref{sec:variational_approach} focuses on computing Problem~\eqref{eq:research_question_upper_bound} derived from a family of product-functions, which is particularly suited to first-order moment assumptions. We formally compare these bounds to existing results. In Section~\ref{sec:Hoeffding}, we study in depth the variational and separable approaches associated with Hoeffding's inequality. Furthermore, we propose a methodology for reconstructing distributions that match the separable or variational optimization problem in the worst-case scenario. Section~\ref{sec:SoS} considers a family of polynomial upper bounds and approximate the resulting upper optimization Problem~\eqref{eq:research_question_upper_bound} using sum-of-square formulations. Thereby, we derive numerical evaluations of Bennett's and Bernstein's inequality. Finally, we introduce a feature-based framework that expands the scope of upper bounds families, providing a comprehensive approach to Hoeffding's inequality that includes both the polynomial and variational methodologies.

\paragraph{Codes.} All codes are provided at \url{https://github.com/CMoucer/ConcentrationInequalities}. We use standard solvers SCS~\citep{o2016conic} and MOSEK~\citep{mosek}.

\section{A variational approach based on product-functions.}
\label{sec:variational_approach}

Classical functions~$F(\cdot)$ usually correspond to the probability tails of a sum of independent random variables, specifically, $F(X) = \mathbf{1}_{\sum_{i=1}^n X_i \in S}$ representing $\Prob(\sum_{i=1}^n X_i \in S)$. This applies, among others, to Hoeffding's, Bennett's, and Bernstein's inequalities (see, e.g.,~\citep{2013Boucheron, 2024Bach, vershynin_2018} and references therein). Their proofs typically rely on the Cramér-Chernoff technique, which essentially combines the exponential Chernoff's inequality with the independence of random variables. Specifically, they use moment-generating functions as follows:
\begin{equation*}
    \Prob\left(\sum_{i=1}^n (X_i - \mu_i) \geqslant nt\right) \leqslant \inf_{\lambda \geqslant 0} e^{-\lambda nt}\E[e^{\lambda \sum_{i=1}^n (X_i - \mu_i)}] = \inf_{\lambda \geqslant 0} \prod_{i=1}^ne^{-\lambda(\mu_i + t)}\E[e^{\lambda X_i}],
\end{equation*}
and then optimize over $\lambda\geqslant 0$ for obtaining the smallest possible valid upper bound within this family. A natural, but richer, family of inequalities for obtaining concentration bounds involves constructing upper bounding functions as products of univariate functions, which are classical in probabilistic variational inference:
\begin{equation}
    \label{eq:upper_functions_product}
    \mathcal{U}= \left\{U: \X^n \mapsto \R^+ {\rm \ such \ that \ } \forall x \in \X^n, U(x) = \prod_{i=1}^n u_i(x_i) {\rm \ and \ } \forall i, u_i : \X \mapsto \R^+ \right\}.
\end{equation}
In the Cramér-Chernoff technique, the functions $u_i$ correspond to moment-generating functions with $\forall x_i \in \X, \ u_i(x_i) = e^{\lambda (x_i - \mu_i - t)}$ and serve as natural upper bounds to the indicator function $\mathbf{1}_{\sum_{i=1}(x_i - \mu_i) \geqslant n t} \leqslant \prod_{i=1}^n e^{\lambda (x_i - \mu_i - t)}$.

Throughout this section, we assume the existence of a product-function $U \in \mathcal{U}$~\eqref{eq:upper_functions_product} such that $\forall x \in \X^n, F(x) \leqslant U(x)$, from which we define two strategies. First, we introduce a separable approach that generalizes classical probability proofs. Second, we optimize over the family of product-functions~\eqref{eq:upper_functions_product} and show how it formulates as a convex optimization problem. Depending on the moments and product-functions under consideration, we demonstrate that these approaches yield tractable upper bounds.

\subsection{The separable approach.}

Let $U \in \mathcal{U}$ be a product-function verifying $\forall x \in \X^n, \ F(x) \leqslant U(x)$. This section studies the properties of the optimization problem:
\begin{equation*}
    \begin{aligned}
        \rho_n^{U} = \sup_{\forall i, p_i \in \Px_{\mu_i}(\X)} \int_{\X^n} U(x_1, \ldots x_n) dp_1(x_1) \cdots dp_n(x_n).
    \end{aligned}
\end{equation*}
This problem naturally provides an upper bound to the generalized problem of moments for independent random variables~\eqref{eq:research_question}, specifically $\rho_n \leqslant \rho_n^{U}$. By definition, the function $U \in \mathcal{U}$ has a product structure $\forall x \in \X^n, U(x) = \prod_{i=1}^n u_i(x_i)$. This allows the separation of the integral over $\X^n$ into $n$ integrals over $\X$, thus decoupling the optimization problem into $n$ independent optimization problems over $p_i \in \Px(\X)$. In other words, the family $\mathcal{U}$ aligns with the structure imposed by independence:
 \begin{equation*}
     \begin{aligned}
          \rho_n^U = \sup_{\forall i, p_i \in \Px_{\mu_i}(\X)} \int_{\X} \cdots \int_{\X} \prod_{i=1}^n u_i(x_i) dp_1(x_1) \cdots dp_n(x_n) = \prod_{i=1}^n \sup_{p_i \in \Px_{\mu_i}(\X)} \int_{\X}u_i(x_i) dp_i(x_i).
     \end{aligned}
 \end{equation*}
For univariate distributions with $\mu_i$ in the interior of $\X$, strong duality holds. Then,
\begin{equation}
\label{eq:exp_dual_reformulation}
\begin{aligned}
     \rho_n^U &=  \prod_{i=1}^n \inf_{\alpha_i \in \R, \beta_i \in \R^m}\{\alpha_i + \beta_i^\top \mu_i \} {\rm \ such \ that \ } \forall \ x_i \in \X, u_i(x_i) \leqslant \alpha_i + \beta_i^\top g_i(x_i), \\
    & = \prod_{i=1}^n \inf_{\alpha_i \in \R, \beta_i \in \R^m} \sup_{x_i \in \X} \{ u_i(x_i) - \beta_i^\top (g_i(x_i) - \mu_i)\}.
    \end{aligned}
\end{equation}
Problem~\eqref{eq:exp_dual_reformulation} is finite-dimensional and convex, as it is the pointwise supremum of affine functions~\cite[Section 3.2.3]{Boyd2004}. In the case of i.i.d. random variables, the problem simplifies significantly to a single optimization problem with $ \rho_n^U = \left( \rho_1^{{\rm exp}}\right)^n$. However, computing $\sup_{x_i \in \X}\{u_i(x_i) - \beta_i^\top(g_i(x_i) - \mu_i)\}$ often remains numerically intractable. Under strong assumptions, such as the finiteness of the support or the convexity of the objective function on a compact set, it reduces to a finite number of constraints. Proposition~\ref{prop:easy_exponential_cases} outlines a useful tractable setting that will be used in Section~\ref{sec:Hoeffding} for improving Hoeffding's inequality.

\begin{proposition}
\label{prop:easy_exponential_cases}
    Let $\X$ be compact and $x_i \mapsto u_i(x_i) - \beta^\top (g_i(x_i) - \mu_i)$ be convex. 
    Then, Problem~\eqref{eq:exp_dual_reformulation} formulates as a tractable convex optimization problem with a finite number of constraints:
    \begin{align*}
     \rho_n^U &=  \prod_{i=1}^n \inf_{\alpha_i \in \R, \beta_i \in \R^m}\{\alpha_i + \beta_i^\top \mu_i \} {\rm \ such \ that \ } \forall x \in {\rm Extremal}(\X), \ \forall i, \  u_i(x_i) \leqslant \alpha_i + \beta_i^\top g_i(x_i).
    \end{align*}
\end{proposition}
\textbf{Proof.} Under the assumptions of Proposition~\ref{prop:easy_exponential_cases}, the maximization of the convex function $u_i(x_i) - \beta_i^\top (g_i(x_i) - \mu_i)$ over the compact set $\X$ is achieved at extremal points of $\X$~\cite[Section 3.2.3]{Boyd2004}. \hfill \Halmos

This technique faces two major challenges: first, constructing a valid upper bound $U \in \mathcal{U}$ can be difficult; second, even with a suitable upper bound, the resulting optimization Problem~\eqref{eq:exp_dual_reformulation} may be numerically intractable. The next section focuses on more accurate approximations to the generalized problem of moments by optimizing over the family of upper bounds (instead of keeping one such upper bound fixed).

\subsection{Optimizing over $\mathcal{U}$.}
This section explores Problem~\eqref{eq:research_question_upper_bound} for the class of product-functions~\eqref{eq:upper_functions_product}. Combined with the dual formulation~\eqref{eq:exp_dual_reformulation}, we define: 
\begin{equation}
\label{eq:definition_rho_var_problem}
\begin{aligned}
     \rho_n^{{\rm var}} = \inf_{u_i \geqslant 0} \inf_{\alpha \in \R^n, \beta\in \R^{n \times m}}  \ \prod_{i=1}^n (\alpha_i + \beta_i^\top\mu_i)  {\rm \ such \ that \ } & \forall i, \forall x_i \in \X, u_i(x_i) \leqslant \alpha_i +\beta_i^\top (x_i), \\ &\forall x \in \X^n, F(x_1, \ldots, x_n) \leqslant \prod_{i=1}^n u_i(x_i).
\end{aligned}
\end{equation}
Optimizing with respect to $u_1, \ldots, u_n \geqslant 0$, it holds that:
\begin{equation}
\label{eq:raw_variational_multivariate_general}
    \begin{aligned}
         \rho_n^{{\rm var}} = \inf_{\alpha \in \R^n, \beta\in \R^{n \times m}} \ \prod_{i=1}^n (\alpha_i + \beta_i^\top \mu_i) {\rm \ such \ that \ } &\forall x_i \in \X, \forall i, \alpha_i + \beta_i^\top g(x_i) \geqslant 0, \\
    & \forall x \in \X^n, F(x_1, \ldots, x_n) \leqslant \prod_{i=1}^n (\alpha_i + \beta_i^\top g(x_i)).
    \end{aligned}
\end{equation}
As expected, Problem~\eqref{eq:raw_variational_multivariate_general} shows no dependence on the $u_i$'s. In Proposition~\ref{prop:comparison_var_exp}, we formally compare $\rho_n$~\eqref{eq:research_question} to~$\rho_{n}^{\rm var}$~\eqref{eq:raw_variational_multivariate_general}~and to $\rho_{n}^{U}$~\eqref{eq:exp_dual_reformulation} for any function $U \in \mathcal{U}$.

\begin{proposition}
    \label{prop:comparison_var_exp}
    Let $U \in \mathcal{U}$, $\rho_n^U$ be defined in~\eqref{eq:exp_dual_reformulation}, $\rho_n$~in~\eqref{eq:research_question} and $\rho_n^{{\rm var}}$~in~\eqref{eq:raw_variational_multivariate_general}. Then it holds that 
    \begin{equation*}
        \rho_n  \leqslant \rho_{n}^{\rm var}\leqslant  \rho_n^U.
    \end{equation*}
    In addition, the equality $ \rho_n^U =  \rho_n^{{\rm var}}$ holds for optimal values $(u_\star, \alpha_\star, \beta_\star)$ such that $\forall x_i \in \X, u_{i, \star}(x_i) = \alpha_{i,\star} + (\beta_{i,\star})^\top g_i(x_i)$.
    \end{proposition}

Proposition~\ref{prop:comparison_var_exp} provides optimal product-functions $U \in \mathcal{U}$, as affine functions of the moments. This specific structure indicates that moment-generating functions, used in the Cramér-Chernoff method, are not optimal. Computing $(\alpha_{i, \star}, \beta_{i, \star})$ often remains difficult (that is, solving Problem~\eqref{eq:raw_variational_multivariate_general} which is finite-dimensional but nonconvex). Proposition~\ref{prop:var_approach_convexified} ensures a convex reformulation of Problem~\eqref{eq:raw_variational_multivariate_general}.

\begin{proposition}
\label{prop:var_approach_convexified}
Let $ \rho_n^{{\rm var}}$ be defined in~\eqref{eq:definition_rho_var_problem}. Then, it holds that
    \begin{align}
\label{eq:rho_var_convexified}
    \log( \rho_n^{{\rm var}}) = \inf_{\alpha \in \R^n, \beta \in \R^{n \times m}} \sum_{i=1}^n \{\alpha_i + \beta_i^\top \mu_i - 1\} + \sup_{x \in \X^n} \left\{\log(F(x)) - \sum_{i=1}^n \log(\alpha_i + \beta_i^\top g_i(x_i))\right\}.
\end{align}
In addition, if $ x \mapsto \log(F(x)) - \sum_{i=1}^n \log(\alpha_i + \beta_i^\top g_i(x_i))$ is convex and $\X$ is compact, then, 
\begin{equation*}
    \sup_{x \in \X^n} \left\{\log(F(x)) - \sum_{i=1}^n \log(\alpha_i + \beta_i^\top g_i(x_i))\right\} = \sup_{x \in {{\rm Extremal}(\X^n)}} \left\{\log(F(x) - \sum_{i=1}^n \log(\alpha_i + \beta_i^\top g_i(x_i))\right\}.
\end{equation*}
\end{proposition}
\textbf{Proof.} First, let us consider the logarithm $\log(\rho_n^{{\rm var}}) = \inf_{\alpha \in \R^n, \beta \in \R^{n \times m}} \sum_{i=1}^n \log(\alpha_i + \beta_i^\top \mu_i) + \sup_{x \in \X^n} \{\log(F(x) - \sum_{i=1}^n \log(\alpha_i + \beta_i^\top g_i(x_i))\}$.
Noticing that $ \inf_{t \geqslant 0}\{\left(t\alpha_i + t\beta_i^\top \mu_i - 1\right) - \log(t\alpha_i + t\beta_i g(x_i)) \} = \log(\alpha_i + \beta_i^\top \mu_i) - \log(\alpha_i + \beta_ig(x_i))$, we conclude the reformulation in~\eqref{eq:rho_var_convexified}. The second assertion follows by maximization of a convex function over a compact set. \hfill \Halmos

Proposition~\ref{prop:var_approach_convexified} details a set of assumptions on $\X$, $F$ and $g$, under which the constraints $ ``\log(F(x_1, \ldots, x_n) \leqslant \sum_{i=1}^n \log(\alpha_i + \beta_i^\top g_i(x_i))$'' reduces to a finite number of points. For instance, these assumptions are satisfied for finite first-order moments $g_i(x_i) = x_i$ together with log-convex objectives $F$ (such as exponentials or indicator functions $\mathbf{1}_{S}$, with compact sets $S \subset \X^n$, see~\cite[Section 3.5]{Boyd2004}). An alternative convexification proof is achieved via optimal transport in Appendix~\ref{app:optimal_transport}, which also includes a formulation of the gap to the generalized problem of moments.

We have established two approaches for deriving upper bounds to the generalized problem of moments for independent random variables~\eqref{eq:research_question} using a family of product-functions~\eqref{eq:upper_functions_product}. First, we introduced a separable approach~\eqref{eq:exp_dual_reformulation} that formulates as a product of $n$ convex optimization problems. However, constructing product functions may not be straightforward. Then, we formulate a variational Problem~\eqref{eq:raw_variational_multivariate_general} emerging from~\eqref{eq:exp_dual_reformulation} by optimizing with respect to product-functions. It turns out that Problem~\eqref{eq:raw_variational_multivariate_general} benefits from a convex reformulation which does not require a priori upper bounds to $F$ (but constructs such bounds in the process). Without further assumptions on the probability support~$\X$, the objective~$F$ or moments~$g_i$, both approaches are intractable. We will see next how they effectively apply in the context of Hoeffding's inequality.
   
\section{Revisiting Hoeffding's inequality.}
\label{sec:Hoeffding}

 Hoeffding's inequality establishes a subgaussian tail for the sum of independent random variables taking their values in a bounded set with finite means. This section is devoted to refining Hoeffding's inequality, applying the separable and variational frameworks developed in Section~\ref{sec:variational_approach}. First, let us recall Hoeffding's inequality as stated in Theorem~\ref{th:Hoeffding}.
\begin{theorem}{\citep{1963Hoeffding}}
\label{th:Hoeffding}
    Let $X_1, ..., X_n$ be independent random variables taking their values in $[a_1, b_i]$ almost surely. Then, for every $t \geqslant 0$,
    \begin{equation}
    \label{eq:Hoeffding_value}
        \Prob\left(\frac{1}{n}\sum_{i=1}^n(X_n - \E[X_n]) \geqslant t\right) \leqslant \exp\left(-\frac{2n^2t^2}{\sum_{i=1}^n(b_i - a_i)^2}\right) = \rho_n^{{\rm Hoeffding}}.
    \end{equation}
\end{theorem}
 Throughout this section, $X_1, \ldots, X_n$ are independent random variables with mean $\mu_1, \ldots, \mu_n$. Without loss of generality, we assume that they all take their values in $\X = [0, 1]$ (that is, $b_i = 1$ and $a_i = 0$) and thereby, have finite mean $\E[X_i] = \mu_i$ for all $i$. Our goal is to approximate the probability $\Prob(\sum_{i=1}^n X_i \geqslant nt + \sum_{i=1}^n \mu_i)$ for $t\geqslant 0$. In our framework, it translates to functions where for all $x \in [0, 1]^n, F(x) = \mathbf{1}_{\sum_{i=1}^n X_i - \sum_{i=1}^n \mu_i - nt}$ and $\forall i, \forall x_i \in [0, 1], g_i(x_i) = x_i$. 

As a reference, we first consider one random variable, that benefits from an exact analytical bound~\eqref{eq:research_question}. We compare it to the separable technique~\eqref{eq:exp_dual_reformulation} on the moment-generating function. Then, we extend the analysis to $n$ random variables, comparing bounds in the separable and variational approaches to Hoeffding's inequality. More precisely, we examine cases where random variables are i.i.d., and where random variables are divided into two blocks with different means $\mu_1 = \mu_2 = \cdots = \mu_m$ and $\mu_{m+1} = \cdots = \mu_n$ (with $1 \leqslant m \leqslant n-1)$. For the case $\mu_1 = \mu_2 = \cdots = \mu_n$, our results asymptotically correspond to large deviations. Finally, we propose a methodology to reconstruct a distribution in the worst-case scenario.

\subsection{One random variable: comparison of the exact and exponential bounds}
\label{sec:hoeffding_one_variable}

Computing optimal bounds for univariate random variables has been extensively studied in past years and is encompassed in the multivariate analyses proposed by~\cite{1962_Isii, 2000_Bertsimas_Popescu_Sethuraman, 2007_Vandenberghe_Boyd_Comanor}. Let us compute the exact closed-form solution~\eqref{eq:research_question} and the separable scenario for moment-generating functions~\eqref{eq:exp_dual_reformulation}.

\paragraph{Exact optimization problem.} Let $X_1$ be a random variable in $\X = [0, 1]$, with $\E[X_1] = \mu$. Recall the exact optimization problem~\eqref{eq:research_question} for every $t \geqslant 0$,
\begin{equation}
\label{eq:exact_univariate}
\begin{aligned}
  \rho_1(t) &= \sup_{p_1 \in \Px([0, 1])} \int_0^1 \mathbf{1}_{x_1 \geqslant \mu + t} dp_1(x_1) \ {\rm \ such \ that \ } \int_{0}^1 x_1 dp_1(x_1) = \mu, \\
    &= \inf_{\alpha, \beta \in \R} \alpha + \beta \mu \ {\rm \ such \ that \ } \forall \ x_1 \in [0, 1], \mathbf{1}_{x_1 \geqslant \mu + t} \leqslant \alpha + \beta x_1.
\end{aligned}
\end{equation}
Problem~\ref{eq:exact_univariate} defines a function $\rho_1(\cdot)$ as a solution to a linear program for every $t \geqslant 0$ and verifies $\rho_1(t) = \Prob\left(X_1 \geqslant t + \mu\right)$. As stated by~\citet[Theorem 2.2]{2005Bertsimas_popescu}, strong duality holds for $\mu \in ]0, 1[$. It turns out that Problem~\eqref{eq:exact_univariate} is a particular case of both the generalized moment problem~\eqref{eq:multivariate_upper_bound} for univariate distributions and of the generalized problem of moments for independent random variables~\eqref{eq:independent_upper_bound_problem}. It admits a closed-form solution, detailed in Proposition~\ref{prop:closed_form_univariate}.

\begin{proposition}
\label{prop:closed_form_univariate}
    Let $0 \leqslant t \leqslant 1 - \mu$. Then, $\rho_1(t)$ as defined in~\eqref{eq:exact_univariate} verifies:
    \begin{equation}
    \label{eq:exact_univariate_value}
        \rho_1(t) = \frac{\mu}{\mu + t}.
    \end{equation}
\end{proposition}
\textbf{Proof.} Functions $x_1 \mapsto -(\alpha + \beta x_1)$ is convex on $[0, 1]$ . Thus the constraint in~\eqref{eq:exact_univariate} can be reduced to two constraints : $\alpha \geqslant 0$ and $\beta (\mu + t) \geqslant 1$. It follows that $\alpha = 0$ and $\beta = \frac{1}{\mu + t}$. \hfill \Halmos

\paragraph{Separable approach.} Now, let us compute the bound defined in~\eqref{eq:exp_dual_reformulation} considering moment-generating functions $u_\lambda(x) = e^{\lambda(x - (\mu + t))}$. By construction, $\forall x \in \R,  \mathbf{1}_{x \geqslant \mu + t} \leqslant e^{\lambda(x - (\mu + t))} $. We define the family of upper bounds $ \rho_1^{{\rm exp}}$. For every $\lambda \in \R$  and for every $t\geqslant 0$,
\begin{equation}
\label{eq:univariate_exp_lambda}
     \rho_1^{{\rm exp}}(\lambda, t) = \inf_{\alpha, \beta} \alpha + \beta \mu, {\rm \ such \ that \ } \forall x_1 \in [0, 1], e^{\lambda(x_1 - (\mu + t))} \leqslant \alpha + \beta x_1.
\end{equation}
Problem~\eqref{eq:univariate_exp_lambda} is a convex optimization problem, that is well-defined for every $t \geqslant 0$ and $\lambda \in \R$. For every $t \geqslant 0$, it admits an optimal moment-generating function analytically given in Proposition~\ref{prop:univariate_exponential}.
\begin{proposition}
\label{prop:univariate_exponential}
    Let $t \geqslant 0$, and $ \rho_1^{{\rm exp}}(\lambda)$ be defined in~\eqref{eq:univariate_exp_lambda}. Then, for every $t \geqslant 0$,
    \begin{equation}
    \label{eq:exp_n1}
        \rho_{1}^{ {\rm exp}}(t) = \left(\frac{\mu}{\nu}\right)^\nu \left(\frac{1 - \mu}{1 - \nu}\right)^{1 - \nu},
    \end{equation}
    where $\nu = \mu + t$, and $\rho_{1, \star}^{\rm exp}$ is the optimal value when optimizing~\eqref{eq:univariate_exp_lambda} with respect to $\lambda$.
\end{proposition}
\textbf{Proof}. See Appendix~\ref{app:proof_univariate_exonential}. \hfill \Halmos

Proposition~\ref{prop:univariate_exponential} yields exactly the Chernoff-Bound for a Bernoulli random variable on the support $\{0, 1\}$ with mean $\mu$~\cite[Section 2.2]{2013Boucheron}. This result was originally stated in the seminal work of \cite{1963Hoeffding}, but our approach ensures tightness for this family of moment-generating functions. It can also be directly derived from Kullback's inequality, as detailed in Appendix~\ref{app:proof_exp_univariate_alternative}.

\begin{figure}[h]
     \begin{subfigure}{0.47\textwidth}
         \centering
         \includegraphics[height=150pt]{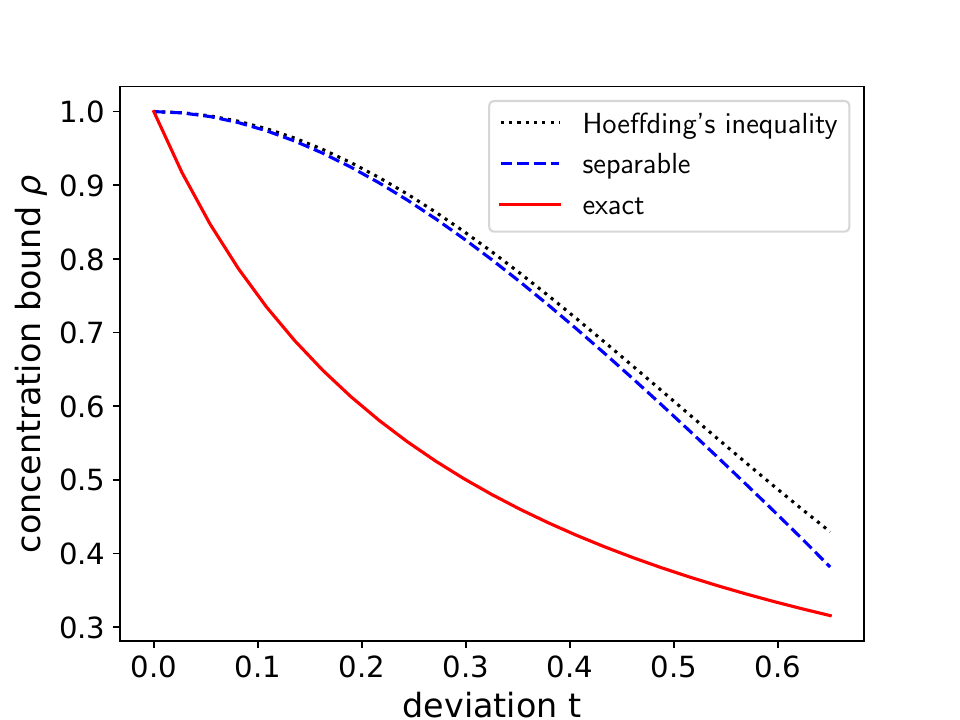}
         \label{fig:comparison_mu}
     \end{subfigure}
     \hfill
     \begin{subfigure}{0.47\textwidth}
         \centering
         \includegraphics[height=150pt]{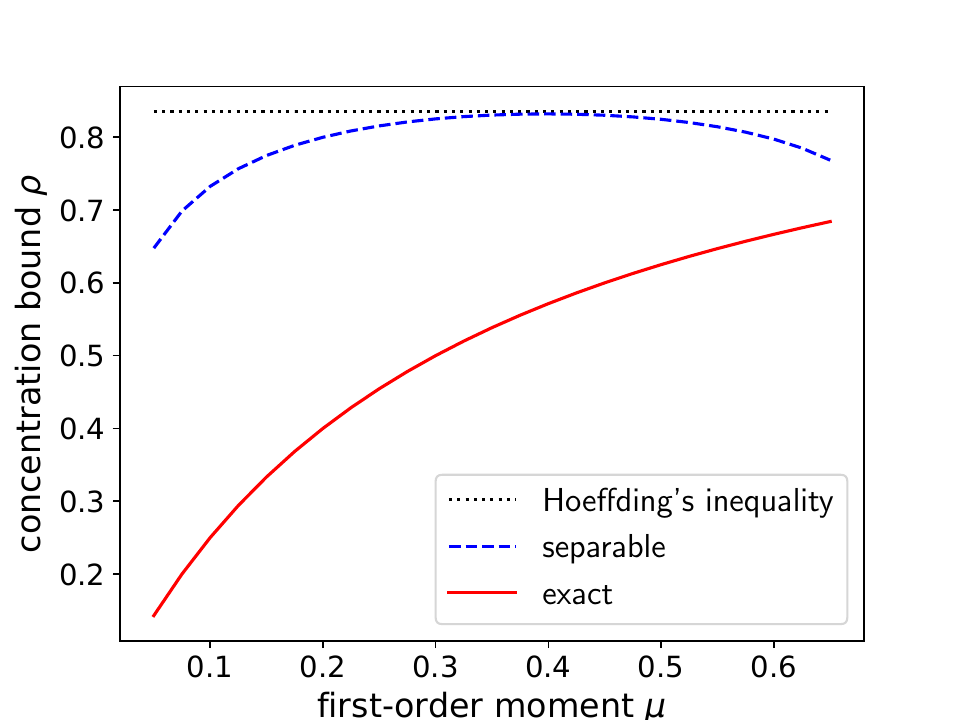}
         \label{fig:comparison_t}
     \end{subfigure}
     \caption{Comparison of the bound in the exact~\eqref{eq:exact_univariate_value} and separable~\eqref{eq:exp_n1} approaches to Hoeffding's inequality~\eqref{eq:Hoeffding_value}. On the left, bounds are plot as a function of $t$ with $\mu=0.3$, and on the right as a function of $\mu$ with $t = 0.3$.}
     \label{fig:comparison_exact_hoeffding_n1}
\end{figure}
In Figure~\ref{fig:comparison_exact_hoeffding_n1}, we observe that the exact bound $\rho_1$~\eqref{eq:exact_univariate} significantly improves upon Hoeffding's bound, whereas the bound $\rho_1^{{\rm exp}}$~\eqref{eq:univariate_exp_lambda} in the separable approach only shows improvement for large values of $t$. Both approaches benefit from a dependence in the first-order moment $\mu$, that does not appear in Hoeffding's inequality~\eqref{eq:Hoeffding_value}. The next section extends beyond the univariate case.

\subsection{Generalization to $n$ independent random variables.}

We consider $X_1, \ldots, X_n$ independent random variables with finite means $\mu_i$. We first examine the case of i.i.d. random variables, where $\mu_1 = \cdots = \mu_n$. The bound obtained in the variational approach~\eqref{eq:raw_variational_multivariate_general} asymptotically matches the large deviations, and cannot therefore be much improved for a large number of variables. Next, we consider the case of independent random variables with different means. Specifically, we formulate a tractable upper bound using the separable treatment and discuss the computational limits of the variational approach for large $n$ due to an exponential number of constraints. In the special case of two blocks of variables with different means, the number of constraints involved is~$O(n^2)$.

\subsubsection{Independent and identically distributed random variables (equal means).}

This section focuses on i.i.d. random variables, that is with $\mu_1 = \mu_2 = \cdots = \mu_n$. Thanks to symmetry properties, both the separable and variational approaches yield tractable solutions for any number of variables~$n$.

\paragraph{Separable approach.} Let us consider the moment generating function of $X_1 + \cdots + X_n$: $\forall \lambda \in \R, \forall x \in [0, 1]^n, u_\lambda(x) = e^{\lambda(\sum_{i=1}^n \{x_i - \mu - t\})} = \prod_{i=1}^ne^{\lambda(x_i - \mu - t)}$. This is a product-function verifying for all $x \in [0, 1]^n, F(x) = \mathbf{1}_{\sum_{i=1}^n x_i \geqslant n(t + \mu)} \leqslant u_\lambda(x)$. We thus define the separable Problem~\eqref{eq:exp_dual_reformulation} for $\lambda \in \R$ and $t \geqslant 0$:
\begin{equation}
\label{eq:exp_multivariate_iid}
     \rho_n^{{\rm exp}}(\lambda, t) = \prod_{i=1}^n \inf_{\alpha_i \in \R, \beta_i \in \R}(\alpha_i + \beta_i \mu) {\rm \ such \ that \ } \forall x_i \in [0, 1]^n, e^{\lambda(x_i - \mu - t)} \leqslant \alpha_i + \beta_i \mu.
\end{equation}
Optimizing with respect to $\lambda$, Proposition~\ref{prop:n_exp_value} provides a closed-form as a function of $t$.
\begin{proposition}
\label{prop:n_exp_value}
    Let $\mu_1 = \cdots = \mu_n$ and let $\rho_n^{\rm exp}$ be defined in~\eqref{eq:exp_multivariate_iid}. Then, it holds for all $t \geqslant 0$:
    \begin{equation}
    \label{eq:exponential_multivariate_iid}
        \rho_{n, \star}^{\rm exp}(t) = \inf_{\lambda\in \R} \rho_n^{\rm exp}(\lambda, t) = \left(\frac{\mu}{\nu}\right)^{n\nu} \left(\frac{1 - \mu}{1 - \nu}\right)^{n(1 - \nu)} = (\rho_{1, \star}^{\rm exp}(t))^n.
    \end{equation}
\end{proposition}
\textbf{Proof}. Since $\mu_1 = \cdots = \mu_n$, Problem~\eqref{eq:exp_multivariate_iid} benefits from a symmetry property and simplifies into  $ \forall t \geqslant 0, \ \forall \lambda \in \R, \  \rho_n^{{\rm exp}}(\lambda, t) = \prod_{i=1}^n \inf_{\alpha, \beta \in \R}(\alpha + \beta \mu), {\rm \ s.t. \ } \forall x_i \in [0, 1]^n, e^{\lambda(x_i - \mu - t)} \leqslant \alpha + \beta \mu$, that is $\rho_n^{{\rm exp}}(\lambda, t) = (\rho_1^{{\rm exp}}(\lambda))^n$. We then optimize over $\lambda$ as in Proposition~\ref{prop:univariate_exponential}. \hfill \Halmos

The symmetry properties induced by the i.i.d. assumption allow simplifying~\eqref{eq:exp_multivariate_iid} into a single univariate convex optimization problem, which can be solved efficiently. Again, this is exactly the Chernoff bound for $n$ i.i.d. Bernoulli variables taking their values in $\{0, 1\}$ with mean $n$. We conclude that the separable approach improves Hoeffding's inequality for large deviations $t$ (as for a unique univariate variable). Let us now explore how the variational approach might offer further improvements.

\paragraph{Variational approach.} In the context of Hoeffding's inequality and given symmetry properties induced by $\mu_1 = \cdots = \mu_n = \mu$, the bound $\rho_n^{{\rm var}}$ defined in~\eqref{eq:raw_variational_multivariate_general} takes the form for all $t \geqslant 0$:
\begin{equation*}
\begin{aligned}
     \rho_n^{{\rm var}}(t) = \inf_{\alpha \in \R, \beta\in \R}  \ \prod_{i=1}^n (\alpha + \beta^\top \mu) \
    {\rm \ such \ that \ } & \forall x \in [0, 1]^n, \ \mathbf{1}_{x_1 + \cdots + x_n \geqslant n(\mu + t)} \leqslant \prod_{i=1}^n (\alpha +\beta x_i), \\
    & \forall x_i \in [0, 1], \ \alpha + \beta x_i \geqslant 0.
\end{aligned}
\end{equation*}
This problem has a convex reformulation with a finite number of constraints, as proven in Proposition~\ref{prop:var_approach_convexified},
\begin{equation}
\label{eq:rho_var_conv_iid}
\begin{aligned}
    \log{ \rho_n^{{\rm exp}}}(t) = \inf_{\alpha, \beta, t \geqslant 0} \  n(\alpha + \beta \mu - 1)
    {\rm \ such \ that \ } & \ -\sum_{i=1}^n \log(\alpha + \beta x_i) \leqslant 0, \ x \in {\rm extremal}(\bar{\X}_n),
\end{aligned}
\end{equation}
where $\bar{\X}_n = \left\{ (x_1, \ldots, x_n) \in [0, 1]^n,  x_1 + \cdots + x_n \geqslant n(\mu + t) \right\}$ and where the constraint ``$\forall x \in [0, 1], \alpha + \beta x \geqslant 0$'' is implied by the logarithm. The set $\bar{\X}_n \subset [0, 1]^n$ is compact and symmetric in $(x_1, \ldots, x_n)$. At first sight, the set ${\rm extremal}(\X)_n$ appears to grow exponentially with $n$. However, due to symmetry properties and the structure of the constraints, it reduces to $O(n)$ constraints (see Appendix~\ref{app:extremal_points_equal_means},  via computation of extremal points). Therefore, Problem~\eqref{eq:rho_var_conv_iid} can be efficiently addressed using standard solvers for convex optimization. In addition, it is possible to derive closed-form solutions for $(\alpha_\star, \beta_\star)$, by enumerating all extremal points given $\mu + t$ and solving the KKT condition for Problem~\eqref{eq:rho_var_conv_iid}. We provide an example for $n=2$ in Appendix~\ref{app:closed_form_reformulation_n2}.

\begin{figure}[h]
     \centering
     \begin{subfigure}[t]{0.32\textwidth}
         \centering
         \includegraphics[height=120pt]{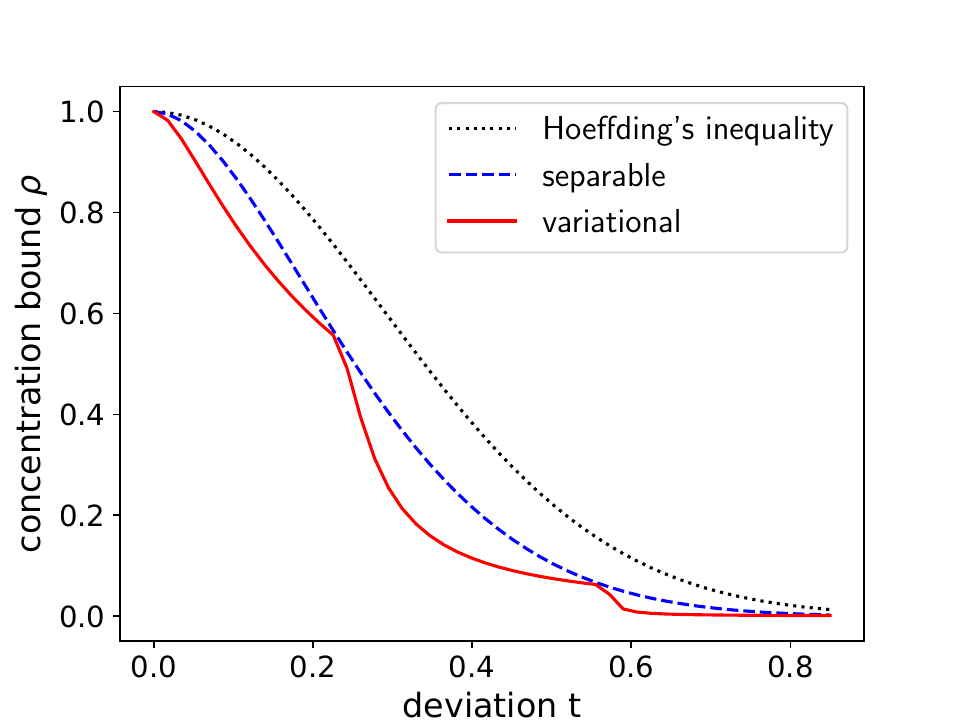}
         \caption{$\mu = 0.1$ and $n = 3$.}
     \end{subfigure}
     \hfill
     \begin{subfigure}[t]{0.32\textwidth}
         \centering
         \includegraphics[height=120pt]{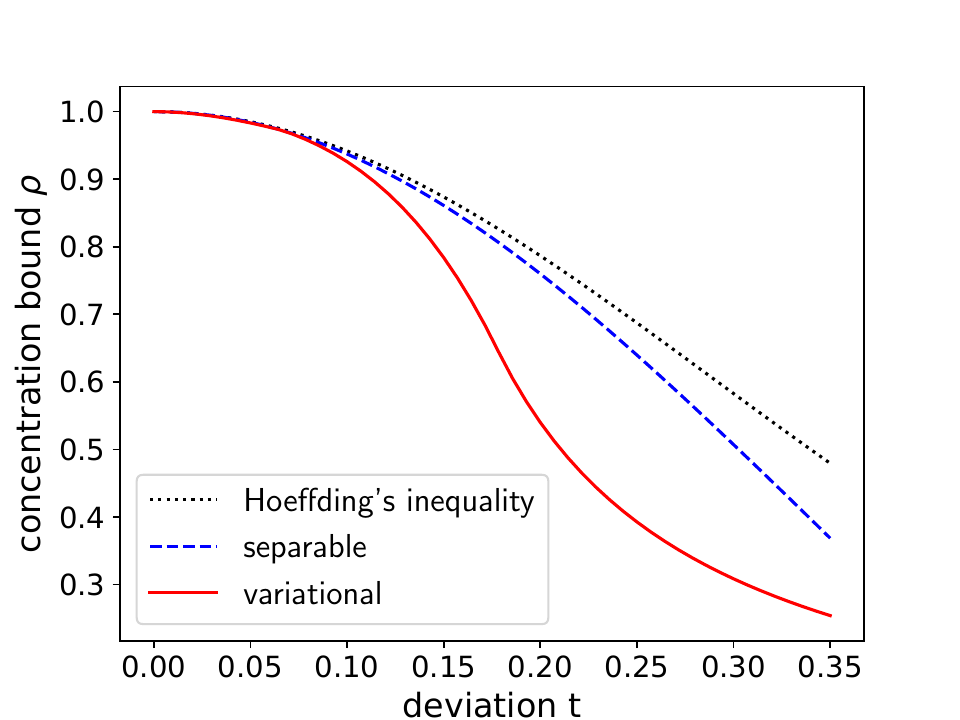}
         \caption{$\mu = 0.6$ and $n=3$.}
     \end{subfigure}
     \hfill
     \begin{subfigure}[t]{0.32\textwidth}
         \centering
         \includegraphics[height=120pt]{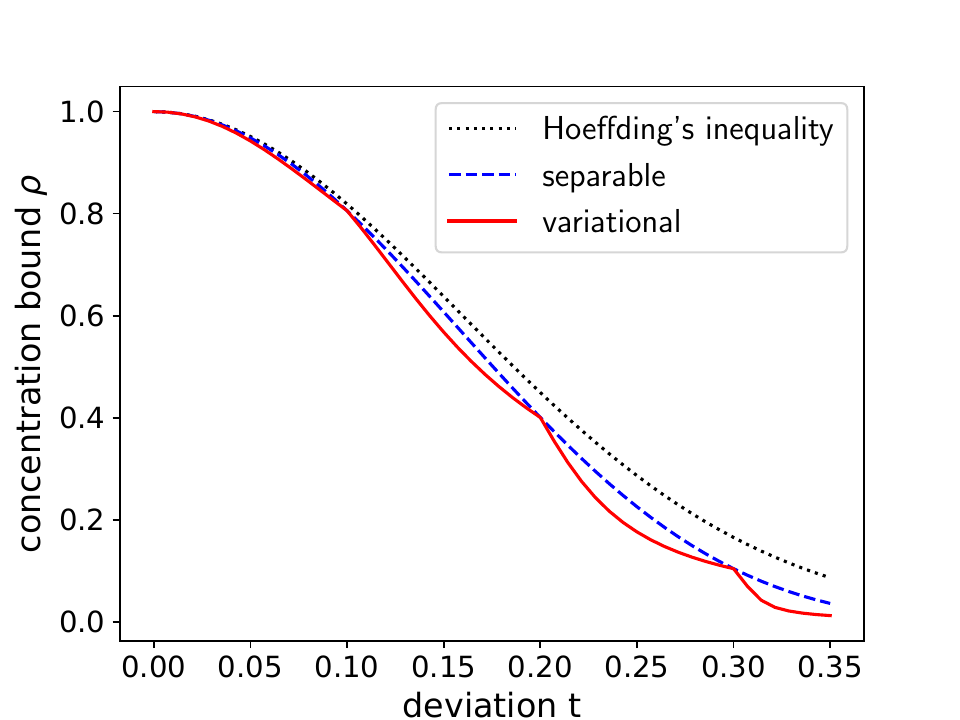}
         \caption{\normalsize{$\mu = 0.6$ and $n = 10$.}}
     \end{subfigure}
     \caption{\normalsize{Comparison of bounds derived in the separable~\eqref{eq:exp_multivariate_iid} and variational approaches~\eqref{eq:rho_var_conv_iid} to Hoeffding's inequality~\eqref{eq:Hoeffding_value}, as a function of the deviation $t$.}}
     \label{fig:comparison_cvx_relaxation_generalization}
\end{figure}

Figure~\ref{fig:comparison_cvx_relaxation_generalization} illustrates the comparison between the separable and variational approaches to the Hoeffding's bound for different numbers of i.i.d. random variables. It appears that $\rho_{n,\star}^{\rm exp}$~\eqref{eq:exp_multivariate_iid} closely tracks Hoeffding's bound. In contrast, $\rho_n^{\rm var}$~\eqref{eq:rho_var_conv_iid} provides significant numerical improvements over Hoeffding's inequality when a small number of random random variables are in play. Furthermore, as the number of variables increases, the variational approach asymptotically matches the separable treatment, as expected from the large deviations theory.

In probability theory, the study of the asymptotic behavior of tails of random variables is known as the large deviations theory, introduced by~\cite{1989_Varadhan}. In particular, the large deviation principle provides a guarantee on rare events, as outlined in Theorem~\ref{th:large_deviations_cramer} for the sum of i.i.d. random variables.

\begin{theorem}[\cite{1938_Cramer}: Large Deviations]
\label{th:large_deviations_cramer}
    Let $X_1, \ldots, X_n$ be i.i.d. random variables with finite moment-generating functions, and let $\bar{X}_n = \frac{1}{n}\sum_{i=1}^n X_i$. Then, for all $x \in \R$,
    \begin{equation*}
        \lim_{n \to \infty} \frac{1}{n}\log(\Prob(\bar{X}_n \geqslant x)) = - \Gamma^\star(x),
    \end{equation*}
    where $\Gamma^\star(x) = \sup_{t \geqslant 0}(tx - \Gamma(t))$ and $\Gamma(t) = \log(\E[\exp(tX_1)])$.
\end{theorem}

The moment-generating function of a univariate random variable formulates as an optimization problem as in~\eqref{eq:univariate_exp_lambda}. Corollary~\ref{coro:large_deviations} provides the large deviations asymptotic for i.i.d. random variables.

\begin{coro}
\label{coro:large_deviations}
    Let $X_1, \ldots, X_n$ be i.i.d. random variables with mean $\E[X_1] = \mu$, and taking their value in $[0, 1]$ almost surely, and let $\bar{X}_n = \frac{1}{n}\sum_{i=1}^n X_i$. It holds that, for all $t \geqslant 0$,
    \begin{equation*}
        \lim_{n \mapsto \infty} \frac{1}{n}\log(\Prob(\bar{X}_n \geqslant \mu + t)) = - (\mu + t)\log\left(\frac{\mu}{\mu + t}\right) - (1 - (\mu + t)) \log\left( \frac{1 - \mu}{1 - (\mu + t)}\right).
    \end{equation*}
\end{coro}
\textbf{Proof.} From the separable approach, $\Gamma(t) = \sup_{p \in \mathcal{P}(\X)} \int_{0}^1 e^{tx}dp(x) \ {\rm such \ that \ } \int_0^1xdp(x) = \mu$. Strong duality holds and $\Gamma(t) = \inf_{\alpha, \beta \in \R} \alpha + \beta \mu, {\rm \ such \ that  \ } \forall x \in [0, 1], \ e^{tx}\leqslant \alpha + \beta x$. Thus, $\Gamma(t) = 1 + \mu(e^t - 1)$. The desired statement is obtained by computing the Fenchel conjugate. \hfill \Halmos

Corollary~\ref{coro:large_deviations} demonstrates that the probability of $\bar{X}_n$ deviating from the mean $\mu$ converges exactly to the bound in the separable approach~\eqref{eq:exponential_multivariate_iid}: for all $t \geqslant 0$, $\Prob(X_1 + \cdots + X_n \geqslant n(\mu + t)) \approx ( \rho_1^{{\rm exp}})^n$. Numerical results presented in Figure~\ref{fig:comparison_cvx_relaxation_generalization} are thus consistent with these large deviation estimates for relatively large~$n$.

As a conclusion, when random variables are i.i.d., optimization problems in the separable and variational approaches benefit from tractable formulations. The bound in the variational treatment~\eqref{eq:rho_var_conv_iid} shows significant improvements over Hoeffding's inequality when a small number $n$ of random variables are in play, but suffers an increasing number of constraints. As $n$ increases, the separable approach~\eqref{eq:exp_multivariate_iid} provides a close estimate of the generalized problem of moments at a lower computational cost.

\subsubsection{Two blocks of random variables with different means.}

Hoeffding's inequality, as presented in Theorem~\ref{th:Hoeffding}, applies generally to independent random variables without specific assumptions on their means. However, formulating the separable and variational approaches in a generic setting can be computationally challenging for a large number of variables $n$. To simplify this, we focus on the scenario where the random variables are divided into two blocks with different means.

\paragraph{Separable approach.} Consider the optimization problem~\eqref{eq:exp_dual_reformulation} with different means in the context of Hoeffding's inequality. After solving each subproblem in $(\alpha_i, \beta_i)$ as in Proposition~\ref{prop:univariate_exponential}, the bound takes the form, for any $t \geqslant 0$:
\begin{equation}
\label{eq:exp_multivariate_noniid}
    \rho_{n, \star}^{\rm exp}(t) = \inf_{\lambda \in \R} e^{-n\lambda(\bar{\mu}_n + t)}\prod_{i=1}^n \left(1 + \mu_i(e^\lambda - 1)\right).
\end{equation}
Optimizing over $\lambda \in \R$ for different means $\mu_i$ cannot be achieved in closed-form as in the case of i.i.d. variables. Note that $\log(\rho_{n, \star}^{\rm exp})$ could be computed as the minimum of a convex objective in $\lambda$. We rather explicit an analytical upper bound in Proposition~\ref{prop:exp_upper_diff_mean}. 
\begin{proposition}
    \label{prop:exp_upper_diff_mean}
    Let $\mu_1, \ldots, \mu_n$ be in $]0, 1[$ and $\mu_i \neq \frac{1}{2}$. Then it holds for $t \geqslant 0$ that: 
    \begin{equation*}
        \rho_{n, \star}^{\rm exp}(t) \leqslant \exp\left( \frac{nt^2/2}{\frac{1}{n} \sum_{i=1}^n \log(\frac{\mu_i}{1 - \mu_i})} \right).
    \end{equation*}
    In addition, $\rho_{n, \star}^{\rm exp}(t) \leqslant \rho_n^{\rm Hoeffding}$.
\end{proposition}
\textbf{Proof. } Let us consider $f_i(\lambda) = \log(1 + \mu_i(e^\lambda - 1))$. A quadratic upper bound for $f$ was derived by \cite[Section 2.2]{Jaakkola2000BayesianPE}, such that $\forall \lambda \in \R, \log(1 + \mu_i(e^\lambda - 1)) \leqslant \lambda \mu_i + \frac{\lambda^2}{4}\frac{2\mu_i-1}{\log(\frac{\mu_i}{1 - \mu_i})}$. In addition, $ \lambda \mu_i + \frac{\lambda^2}{4}\frac{2\mu_i-1}{\log(\frac{\mu_i}{1 - \mu_i})}\leqslant \lambda\mu_i + \frac{\lambda^2}{8}$, leading to the final assertion. \hfill \Halmos
\begin{remark}
In the proof for Proposition~\ref{prop:exp_upper_diff_mean}, considering the naive upper function $\log(1 + \mu_i(e^\lambda - 1)) \leqslant \lambda \mu_i + \frac{\lambda^2}{8}$ would have led to Hoeffding's inequality.
\end{remark}

Given random variables with different means, Proposition~\ref{prop:exp_upper_diff_mean} shows a control of $\rho_{n, \star}^{\rm exp}(t)$ by an upper bound depending on $(\mu_i)_{i=1, \ldots, n}$ and improving Hoeffding's inequality. Therefore, it appears to be suited to different means as well as to two blocks of random variables.

\paragraph{Variational approach.} Recall the optimization problem defining the variational approach for random variables in $[0,1]$ having different means $\mu_1, \ldots, \mu_n$. For $t \geqslant 0$, we define:
\begin{equation}
\label{eq:var_diff_mean}
    \log( \rho_n^{{\rm var}}(t)) = \inf_{\alpha \in \R^n, \beta \in \R^n} \sum_{i=1}^n \{\alpha_i + \beta_i \mu_i -1 \}, {\rm \ such \ that \ } -\sum_{i=1}^n \log(\alpha_i + \beta_i x_i) \leqslant 0, \forall x \in {\rm \ extremal}(\bar{\X}_n),
\end{equation}
where $\bar{\X}_n = \{(x_1, \ldots, x_n) \in [0, 1]^n , \sum_{i=1}^n x_i \geqslant nt + \sum_{i=1}^n \mu_i\}$. Without further assumptions on the $\mu_i$'s, this optimization problem may have up to $O(n!)$ constraints. To simplify the computations, we consider a first group of variables $X_1, \ldots X_m$ with mean $\mu_1$, and a second group $X_{m+1}, \ldots, X_n$ with mean $\mu_2$, with $1 \leqslant m \leqslant n$. Then, the number of constraints under consideration can be reduced to $O(n)$, as shown in Lemma~\ref{lem:nbconstraints_two_blocks}.
\begin{lemma}
    \label{lem:nbconstraints_two_blocks}
    Let $ \mu_1 = \cdots = \mu_m$ and $\mu_{m+1} = \cdots = \mu_n$ in~\eqref{eq:var_diff_mean}. Then, the number of constraints in Problem~\eqref{eq:var_diff_mean} reduces to $O(n^2)$.
\end{lemma}
\textbf{Proof. } See Appendix~\ref{app:proof_lemma_constraints_two_blocks}. \hfill \Halmos

\begin{figure}[h!]
     \centering
     \begin{subfigure}[t]{0.32\textwidth}
         \centering
         \includegraphics[height=125pt]{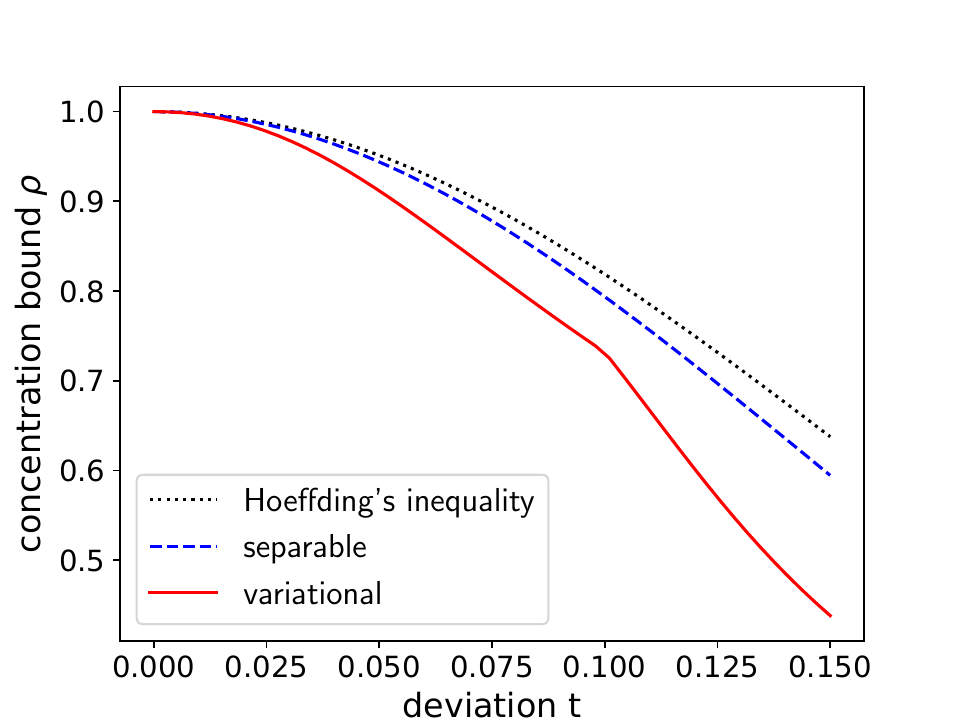}
         \caption{$n=10$, $m=5$, $\mu_1 = 0.2$, $\mu_2 = 0.8$}
     \end{subfigure}
     \hfill
     \begin{subfigure}[t]{0.32\textwidth}
         \centering
         \includegraphics[height=125pt]{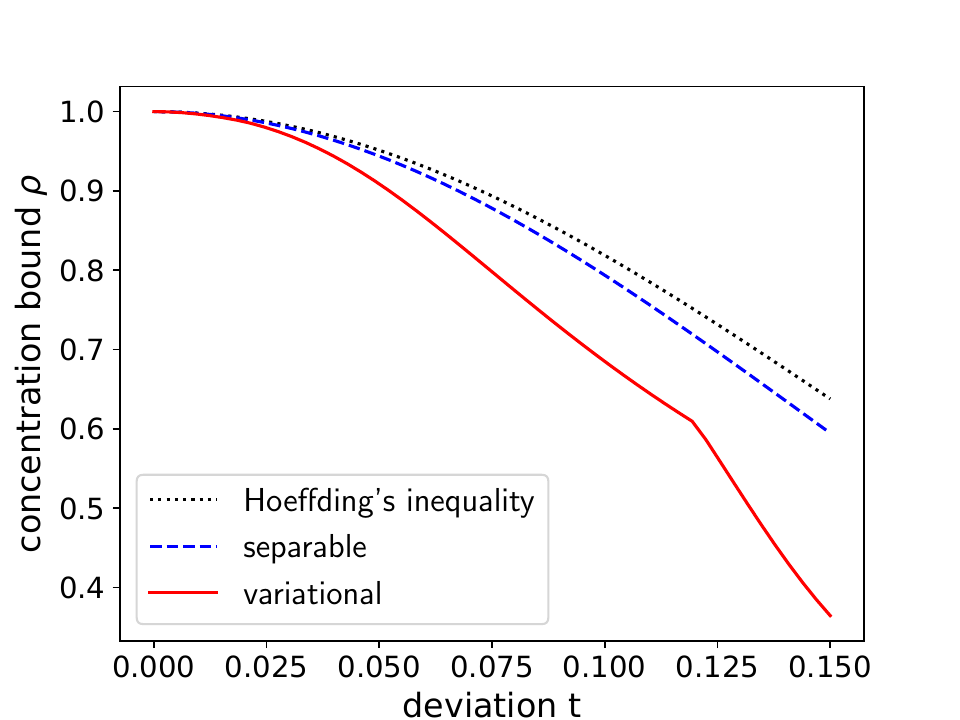}
         \caption{$n=10$, $m=2$, $\mu_1 = 0.2$, $\mu_2 = 0.8$}
     \end{subfigure}
     \hfill 
     \begin{subfigure}[t]{0.32\textwidth}
         \centering
         \includegraphics[height=125pt]{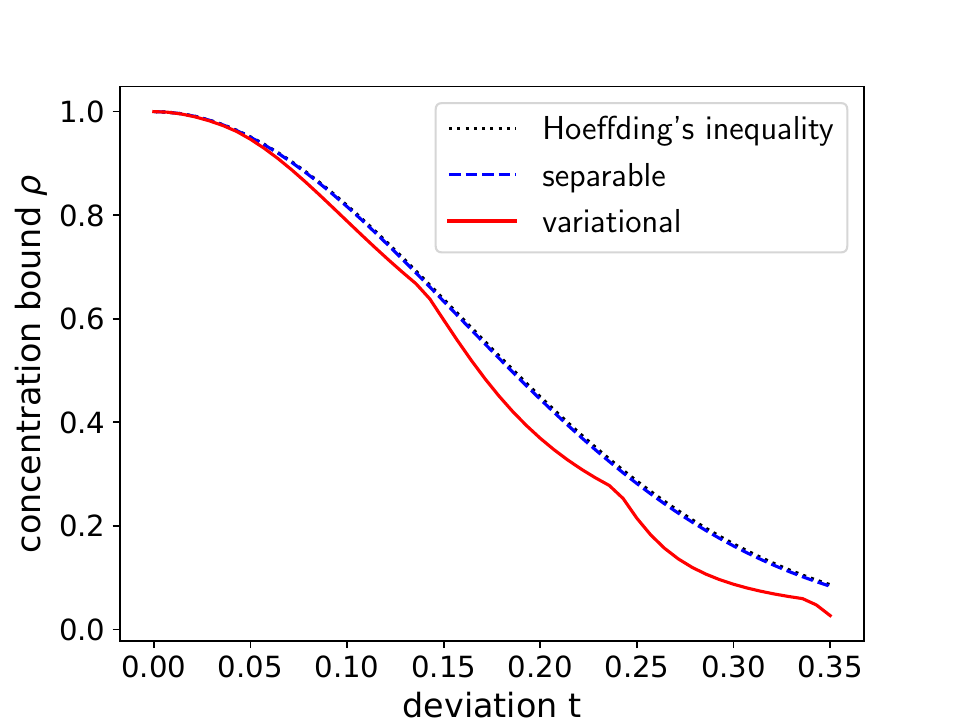}
         \caption{$n = 10$, $m = 2$, $\mu_1 = 0.4$ and $\mu_2 = 0.6$}
     \end{subfigure}
\caption{\normalsize{Comparison of the variational~\eqref{eq:var_diff_mean} and separable~\eqref{eq:exp_multivariate_noniid} approaches to Hoeffding's inequality~\eqref{eq:Hoeffding_value}, in the context of independent random variables divided into two blocks of size $m$ and $n-m$, with mean $\mu_1$ and $\mu_2$,}}
\label{fig:comparison_two_block}
\end{figure}

In Figure~\ref{fig:comparison_two_block}, we observe that the variational approach~\eqref{eq:var_diff_mean} improves largely upon the separable scenario~\eqref{eq:exp_multivariate_noniid}, as soon as the means of the two blocks differs significantly.

Given $n$ random variables divided into two subgroups with different means, we have formulated tractable formulations in the variational~\eqref{eq:var_diff_mean} and separable~\eqref{eq:exp_multivariate_noniid} approaches. On the one hand, we provide an upper bound to the separable approach in Proposition~\ref{prop:exp_upper_diff_mean} using a quadratic upper bound, which still offers improvements over Hoeffding's inequality. On the other hand, the variational approach can be expressed as a manageable convex optimization problem with $O(n)$ constraints, as stated in Lemma~\ref{lem:nbconstraints_two_blocks}. Numerical comparisons show that both approaches behave similarly for a large number of variables $n$, but tend to differ for small $n$. In the following, we explore the reconstruction of extremal distributions involved at the optimum.

\subsection{Reconstructing extremal distributions.}

A bound is said to be tight if a distribution satisfies the bound with equality, as defined by~\citet[Section 2.]{2005Bertsimas_popescu}. Such a distribution is referred to as an \textit{extremal distribution}. We first restrict our attention to one random variable, for which the exact optimization formulation~\eqref{eq:exact_univariate} and the separable approach~\eqref{eq:univariate_exp_lambda} admit analytical solutions. Building on these results, we propose a strategy to construct an extremal distribution for $n$ random variables in the variational~\eqref{eq:var_diff_mean} and separable approaches~\eqref{eq:exp_multivariate_noniid}.

\subsubsection{Dirac distributions for one univariate random variable.}

When it comes to one univariate random variables, we derived in Section~\ref{sec:hoeffding_one_variable} closed form of the exact~\eqref{eq:exact_univariate} and separable~\eqref{eq:univariate_exp_lambda} approaches. In the following, we construct their corresponding extremal distributions.

\paragraph{Exact optimization problem.} Recall the exact optimization problem~\eqref{eq:exact_univariate} for one variable, for all $t \geqslant 0$:
\begin{equation}
\label{eq:exact_dualizing_twice}
\begin{aligned}
    \rho_1(t) &= \sup_{p_1 \in \Px([0, 1])} \int_0^1 \mathbf{1}_{x_1 \geqslant \mu + t} dp_1(x_1) {\rm \ such \ that \ } \int_{0}^1 x_1 dp_1(x_1) = \mu.
\end{aligned}
\end{equation}
Its exact value is provided in Proposition~\ref{prop:closed_form_univariate}, for all $t \geqslant 0, \rho_1(t) = \frac{\mu}{\mu + t}$. Dualizing twice the optimization problem~\eqref{eq:exact_dualizing_twice}, we propose in Proposition~\ref{prop:extremal_distribution_univariate_exact} a strategy for reconstructing an extremal distribution.

\begin{proposition}
    \label{prop:extremal_distribution_univariate_exact}
    Let $t \geqslant 0$. The following distribution is an optimal solution to~\eqref{eq:exact_univariate}:
    \begin{equation}
        \forall x \in [0, 1], p(x) = \frac{t}{\mu + t}\delta_{x = 0} + \frac{\mu}{\mu + t}\delta_{x = \mu + t}.
    \end{equation}
\end{proposition}
\textbf{Proof}. The result is obtained by computing directly $\int_{0}^1\mathbf{1}_{x \geqslant \mu + t} dp(x) = \frac{\mu}{\mu + t}$. We rather propose a constructive approach. First, recall the dual of Problem~\eqref{eq:exact_dualizing_twice}: for all $t \geqslant 0, \rho_1(t) = \inf_{\alpha, \beta \in \R} \alpha + \beta \mu, {\rm \ such \ that \ } \forall x \in [0, 1], \mathbf{1}_{x \geqslant \mu + t} \leqslant \alpha + \beta x$. By convexity of $x \in [0, 1] \mapsto \mathbf{1}_{x \geqslant \mu + t} \leqslant \alpha + \beta x$, the problem reduces for all $t \geqslant 0$ to $\rho_1(t) = \inf_{\alpha, \beta} \alpha + \beta \mu, {\rm \ such \ that \ } \alpha \geqslant 0, 1 \geqslant \alpha + \beta(\mu + t)$. At optimality, $(\alpha_\star, \beta_\star) = (0, \frac{1}{\mu + t})$ meaning that the points $x=0$ and $x=\mu + t$ are active. We compute again the dual, which reformulates as a problem in the probability space by strong duality: for all $t \geqslant 0$ $\rho_1(t) = \sup_{\lambda_1, \lambda_2 \geqslant 0} \lambda_2, {\rm \ such \ that \ } 1 = \lambda_1 + \lambda_2, \mu = \lambda_2(\mu + t)$. The constraint ``$1 = \lambda_1 + \lambda_2$'' corresponds to initial constraint $\int_{0}^1 dp(x) = 1$, and the second one to the first-order moment condition $\int_{0}^1 xdp(x) = \mu$. At optimality, $(\lambda_{1, \star}, \lambda_{2, \star}) = (\frac{t}{\mu + t}, \frac{\mu}{\mu + t})$. We conclude by identification. \hfill \Halmos

\begin{remark}
Other distributions achieve the optimal bound, such as $p(x) = (1 - \mu)\delta_0 + \mu \delta_1$. 
\end{remark}

Proposition~\ref{prop:extremal_distribution_univariate_exact} offers a strategy for reconstructing an extremal distribution, mostly by reducing the constraints involved at the optimum to some active points and by dualizing twice.

\paragraph{Separable approach.} We derive the same technique described in the proof for Proposition~\ref{prop:extremal_distribution_univariate_exact} to the separable approach. Recall the separable approach applied to moment-generating functions~\eqref{eq:univariate_exp_lambda} is defined for $\lambda \in \R$ and $t \in \geqslant 0$,
\begin{equation*}
     \rho_1^{{\rm exp}}(\lambda, t) = \inf_{\alpha, \beta} \alpha + \beta \mu \ {\rm \ such \ that \ } \ \forall \ x_1 \in [0, 1], e^{\lambda(x_1 - (\mu + t))} \leqslant \alpha + \beta x_1.
\end{equation*}
In the proof for Proposition~\ref{prop:univariate_exponential}, we showed that $ \rho_1^{{\rm exp}}(\lambda, t) = e^{-\lambda(\mu + t)}\left((1 + \mu(e^{\lambda} - 1)\right)$. We compute an example of an extremal distribution in Proposition~\ref{prop:extremal_distribution_univariate_exp}, that is independent of $t$ and $\lambda$.

\begin{proposition}
    \label{prop:extremal_distribution_univariate_exp}
    Let $t, \lambda \in \R$. The following distribution is an optimal solution to~\eqref{eq:univariate_exp_lambda}:
    \begin{equation*}
        p(x) = (1 - \mu)\delta_{x=0} + \mu \delta_{x=1}.
    \end{equation*}
\end{proposition}
\textbf{Proof.} Following the same approach as in Proposition~\ref{prop:extremal_distribution_univariate_exact} for determining the extremal distribution of the exact optimization problem. We prove in Appendix~\ref{app:proof_univariate_exonential} that $(\alpha_\star, \beta_\star) = (e^{-\lambda(\mu + t)}, e^{-\lambda(\mu + t)}(e^\lambda - 1))$. After redualizing, it leads to the optimization problem $\sup_{\lambda_1, \lambda_2 \geqslant 0} \lambda_1e^{-\lambda(\mu + t)} + \lambda_2 e^{\lambda(1 - \mu - t)}, {\rm \ such \ that \ } 1 = \lambda_1 + \lambda_2, \mu = \lambda_2$. We conclude by identification. \hfill \Halmos

When studying concentration inequalities applied to a one (univariate) random variable, we derived examples of extremal distributions. This strategy can be decomposed into two steps. First, the convexity properties of the constraints \(\forall x \in [0, 1], \mathbf{1}_{x \geqslant \mu + t} \leqslant \alpha + \beta x\) and \(\forall x \in [0, 1], e^{\lambda (x - \mu - t)} \leqslant \alpha + \beta x\) allow determining active constraints and thus identifying the Dirac delta functions involved at optimality. Second, leveraging strong duality, we dualize the dual, leading to an optimization problem with respect to a simplified space of distributions. We then conclude by identification. In what follows, we show that this technique extends well to $n$ univariate random variables.

\subsubsection{Generalization to $n$ random variables.}

 We now turn to the problem of deriving extremal distribution for $n$ independent random variables. The strategy developed for one random variable extends well to multiple random variables in the separable approach. In the variational approach however, this process often involves computing analytically the solution in the dual before obtaining the extremal distribution.

\paragraph{Separable approach.} The separable approach~\eqref{eq:exp_dual_reformulation} benefits from a decoupling into $n$ independent optimization problems on one random variables. In the context of Hoeffding's inequality, recall its dual formulation~\eqref{eq:exp_multivariate_noniid} below, for $t \geqslant 0$ and $\lambda \in \R$:
\begin{equation*}
     \rho_n^{{\rm exp}}(\lambda, t) = \prod_{i=1}^n \inf_{\alpha_i, \beta_i}(\alpha_i + \beta_i \mu_i) \ {\rm \ such \ that \ } \ \forall x_i \in [0, 1]^n, e^{\lambda(x_i - \mu_i - t)} \leqslant \alpha_i + \beta_i \mu_i.
\end{equation*}
Proposition~\ref{prop:extremal_distribution_univariate_exp} applies on each subproblem $i$, leading to Corollary~\ref{coro:extremal_distribution_exp_multi}.
\begin{coro}
    \label{coro:extremal_distribution_exp_multi}
    The following distribution is an optimal solution to~\eqref{eq:exp_multivariate_noniid} $p(x) = \prod_{i=1}^n p_i(x_i)$, with $\forall x_i \in [0, 1], p_i(x_i) = (1 - \mu_i) \delta_{x_i = 0} + \mu_i \delta_{x_i = 1}.$
\end{coro}

\paragraph{Variational approach.} At first glance, the variational approach reintroduces coupling between variables in the optimization problem. Let us recall its dual formulation~\eqref{eq:var_diff_mean} in the context of Hoeffding's inequality, for $t \geqslant 0$,
\begin{equation}
\label{eq:hoeffding_dual_simplified}
\begin{aligned}
     \rho_n^{{\rm var}}(t) = \inf_{\alpha \in \R^n, \beta \in \R^n} \prod_{i=1}^n (\alpha_i + \beta_i \mu_i) \ {\rm \ such \ that \ } & \forall \ x \in [0, 1]^n, \mathbf{1}_{\sum_{i=1}^n x_i \geqslant \sum_{i=1}^n \mu_i + nt} \leqslant \prod_{i=1}^n (\alpha_i + \beta_i x_i),\\
    &\forall \ i, \forall \ x_i \in [0, 1], 0 \leqslant \alpha_i + \beta_i x_i.
\end{aligned}
\end{equation}
The Lagrangian dual of Problem~\eqref{eq:hoeffding_dual_simplified} cannot be formulated in closed form due to the product form in the constraints. Let us revisit the original optimization problem from which this problem is derived:
\begin{equation}
\label{eq:dual_hoeffding_extremal_distribution}
\begin{aligned}
     \rho_n^{{\rm var}}= \inf_{\forall i, u_i \geqslant 0}\inf_{\alpha \in \R^n, \beta \in \R^n} \prod_{i=1}^n (\alpha_i + \beta_i \mu_i), \
 {\rm \ such \ that \ } & \forall x \in [0, 1]^n, \mathbf{1}_{\sum_{i=1}^n x_i \geqslant \sum_{i=1}^n \mu_i + nt} \leqslant \prod_{i=1}^n u_i(x_i), \\
 &\forall i, \forall x_i \in [0, 1], u_i(x_i) \leqslant \alpha_i + \beta_i x_i.
\end{aligned}
\end{equation}
Using this formulation, Proposition~\ref{prop:extremal_distribution_variational} computes explicitly an extremal distribution.
\begin{proposition}
    \label{prop:extremal_distribution_variational}
    Let $(u_{\star}, \alpha_\star, \beta_\star)$ be a solution to~\eqref{eq:dual_hoeffding_extremal_distribution}. Then,  $\forall i, \forall x_i \in [0, 1], u_{i, \star}(x_i) = \alpha_{i, \star} + \beta_{i, \star}x_i.$
    In addition, an extremal distribution is given by:
    \begin{equation*}
        p(x) = \prod_{i=1}^n \left((1 - \mu_i)\delta_{x_i = 0} + \mu_i \delta_{x_i = 1} \right).
    \end{equation*}
\end{proposition}
\textbf{Proof. } Let $(u_{\star}, \alpha_\star, \beta_\star)$ be optimal solutions in~\eqref{eq:dual_hoeffding_extremal_distribution}. Then, Proposition~\ref{prop:comparison_var_exp} provides the form of the optimal product-function upper bounding $F$, that we recall: $\forall i, \forall x_i \in [0, 1], u_{i, \star}(x_i) = \alpha_{i, \star} + \beta_{i, \star}x_i.$ Thus, considering this specific product-function, we have
\begin{align*}
\rho_{relax}^n &= \prod_{i=1}^n\sup_{p_i \in \mathcal{P}(\X)} \int_{0}^1 u(x_i) dp_i(x_i) \ {\rm such \ that \ } \int_0^1 x_i dp_i(x_i) = \mu_i, \\ 
&= \prod_{i=1}^n \inf_{\lambda_i, \nu_i} (\lambda_i \mu_i + \nu_i) \ {\rm such \ that \ } \alpha_{i, \star} -\nu_i + (\beta_{i, \star} - \lambda_i) x_i \leqslant 0, \forall x_i \in [0, 1],
\end{align*}
where the constraint ``$\alpha_{i, \star} -\nu_i + (\beta_{i, \star} - \lambda_i) x_i \leqslant 0, \forall x_i \in [0, 1]$'' reduces to $\alpha_{i, \star} -\nu_i \leqslant 0$ for $x_i  = 0$ and $\alpha_{i, \star} -\nu_i + \beta_{i, \star} -\lambda_i \leqslant 0$ for $x_i  = 1$. We conclude by identification that $p(x) = (1 - \mu)\delta_0 + \mu \delta_1$. \hfill \Halmos.
 
The extremal distribution derived in Proposition~\ref{prop:extremal_distribution_variational} is exactly equal to the extremal distribution in the separable treatment approach moment-generating functions in Proposition~\ref{coro:extremal_distribution_exp_multi}. Again, it is independent of the deviation $t$. In both cases, it turns out that the dependence in $t$ is only supported in the optimal upper function $U_\star$, either in the moment generating function or in the linear function parametrized by $(\alpha_\star, \beta_\star)$, as detailed in Appendix~\ref{app:closed_form_reformulation_n2} for $n=2$.

In this section, we have thus refined Hoeffding's inequality through two different approaches. When the random variables are i.i.d., the separable approach applied to moment-generating functions aligns with the large deviation principle and slightly improves the traditional Hoeffding's inequality. In contrast, the variational approach yields significantly smaller bounds for a small number of variables but requires $O(n)$ constraints. The case of distinct means $\mu_1, \ldots, \mu_n$ is more advanced and we only explicitly attacked the scenario of two blocks with distinct means. Finally, we proposed a strategy for reconstructing an extremal distribution. As a natural extension, these strategies could potentially be applied to Bennett or Bernstein's inequalities, which assume first and second-order conditions. However, extending to higher-order assumptions reveals challenges where both the separable and variational approaches fail to provide computable solutions. In what follows, we propose a new family of upper bounds adapted to such scenarios.

\section{A polynomial approach based on sum-of-square decomposition.}
\label{sec:SoS}

The separable and variational approaches are constructed from the family of product-functions~\eqref{eq:upper_functions_product} that upper bound $F$ on $\mathcal{X}$. They turn out to be effective for finite first-order moments, but remain computationally out of reach when assuming fixed higher-order moments. Even a simple second-order assumption in Hoeffding's inequality (or an assumption on the variance) results in challenging constraints in the separable approach:
\begin{equation*}
     \rho_n^{{\rm exp}}(\lambda) = \prod_{i=1}^n \inf_{\alpha_i, \beta_i}(\alpha_i + \beta_{i}^{(1)} x_i + \beta_{i}^{(2)}) {\rm \ s.t. \ } \forall x_i \in [0, 1]^n, e^{\lambda(x_i - \mu_i - t)} \leqslant \alpha_i + \beta_{i}^{(1)} x_i + \beta_{i}^{(2)} x_i^2,
\end{equation*}
where $\mu_i^{(1)}$ (resp.~$\mu_i^{(2)}$) represents the first (resp.second) order moment condition. The constraint ``$\forall x_i \in [0, 1]^n, e^{\lambda(x_i - \mu_i - t)} - (\alpha_i + \beta_{i}^{(1)} x_i + \beta_{i}^{(2)} x_i^2)$ admits no closed-form solution. This issue also prevents from computing the large deviations for i.i.d. random variables, which involves computing $\Gamma(t) = \log(\E[\exp(tX_1)])$ as presented in Theorem~\ref{th:large_deviations_cramer}. Similarly, the variational approach involves polynomial constraints with a product structure, making the optimization problem more complex to solve:
\begin{equation*}
    \begin{aligned}
         \rho_n^{{\rm var}} = \inf_{\alpha \in \R^n, \beta\in \R^{n \times m}} \ \prod_{i=1}^n (\alpha_i + \beta_i^{(1)} \mu^{(1)} + \beta_i^{(2)} \mu^{(2)}) &{\rm \ such \ that \ } \forall x_i \in \X, \forall i, \alpha_i + \beta_i^{(1)}x_i + \beta_i^{(2)}x_i^2 \geqslant 0, \\
    & \forall x \in \X^n, F(x_1, \ldots, x_n) \leqslant \prod_{i=1}^n (\alpha_i + \beta_i^{(1)}x_i + \beta_i^{(2)}x_i^2).
    \end{aligned}
\end{equation*}
In this section, we explore polynomial families of upper bounds adapted to such scenarios. First, we propose to analyze Problem~\eqref{eq:research_question_upper_bound} using the family of linear upper bounds, for which we derive closed-form upper bounds. This \textit{linear approach} offers already an improvement to Hoeffding's inequality in comparison with the variational approach for small values of $t$. We then extend this approach to a polynomial family of upper bounds, whose degree equals the number of variables. This results in an optimization problem with an infinite number of polynomial constraints. This so-called \textit{polynomial approach} is closely related to the work of~\cite{2005Bertsimas_popescu}, who introduced a series of SDPs to approximate the generalized problem of moments for multivariate random variables. This approach numerically improves Bernstein's and partially Bennett's inequality. Finally, we introduce a \textit{feature-based approach} generalizing the polynomial approach to a broader family of upper bounds. This allows in particular analyzing Hoeffding's inequality using second-order polynomials.

\subsection{A simple linear upper bound for Hoeffding's concentration.}

Before studying polynomial upper bounds, we start by exploring the simpler family of linear upper bounds and applying it in the context of Hoeffding's inequality. This scenario outlines the key concepts that inspire the polynomial approach and results in an optimization problem that can be solved in closed-form.

 Let $X_1, \ldots, X_n$ be i.i.d. random variables with mean $\E[X_i] = \mu$ and taking their values in $[0, 1]$ and let us introduce the family of linear functions $\mathcal{L} = \{a^\top x + b, a \in \R^n, b \in \R\}$. 
\begin{remark}
    For $n \geqslant 2$, notice that linear functions cannot be formulated as product-functions. Indeed, for two i.i.d. random variables, product-functions takes the form $U(x) = (\alpha + \beta x_1)(\alpha + \beta x_2) = \alpha^2 + \alpha \beta (x_1 + x_2) + \beta^2 x_1 x_2$. 
\end{remark}

The associated optimization problem in the linear approach takes the form, for $t \geqslant 0$:
\begin{equation}
\label{eq:linear_upper_optimization problem}
\begin{aligned}
    \rho_n^{{\rm lin}}(t) = \inf_{a \in \R^n, b \in \R} \sup_{\forall i, p_i \in \Px_{\mu_i}(\X)} & \int_{\X^n} \left(\sum_{i=1}^n a_i x_i + b\right)dp_1(x_1) \cdots dp_n(x_n), \\
    {\rm \ such \ that \ } & \forall x \in \X^n, F(x) \leqslant \sum_{i=1}^n a_ix_i + b.
\end{aligned}
\end{equation}
By construction, $\rho_n \leqslant \rho_n^{{\rm lin}}$. We solve Problem~\eqref{eq:linear_upper_optimization problem} analytically in Proposition~\ref{prop:hoeffding_linear}.

\begin{proposition}
\label{prop:hoeffding_linear}
    Let $X_1, \ldots, X_n$ be i.i.d. random variables taking their values in $[0, 1]$ and with finite mean $\E[X_i] = \mu$. Then, it holds for all $t \in [0, 1 - \mu]$ that:
    \begin{equation*}
        \rho_n^{{\rm lin}}(t) = \frac{\mu}{\mu + t}.
    \end{equation*}
\end{proposition}
\textbf{Proof}. Under the assumptions of Hoeffding's inequality and by symmetry,it holds for $t \geqslant 0$:
\begin{align*}
    \rho_n^{{\rm lin}}(t) &= \inf_{a \in \R, b \in \R} a n \mu + b {\rm \ such \ that \ } \forall x \in [0, 1]^n, \mathbf{1}_{x_1 + \cdots + x_n \geqslant n(\mu + t)} \leqslant a\sum_{i=1}^n x_i + b.
\end{align*}
It follows that $b = 0$ by considering $x = 0$ and minimizing with respect to $b$. By construction, $\sum_{i=1}^n x_i \geqslant n(\mu + t)$ implies $1 \leqslant a\sum_{i=1}^n x_i \leqslant an(\mu + t)$. We conclude that $a = \frac{1}{n(\mu + t)}$ and the desired result. \hfill \Halmos

In Proposition~\ref{prop:hoeffding_linear}, we show that $\rho_n^{{\rm lin}}$ is exactly equal to the exact bound for one univariate random variable as defined in~\eqref{eq:exact_univariate}, that is $\rho_n^{{\rm lin}} = \rho_1$. In addition, $\rho_n^{\rm lin}$ shows a dependence on the mean $\mu$, but no dependence on the number of variables under consideration. 
\begin{figure}[h!]
\centering
\includegraphics[height=150pt]{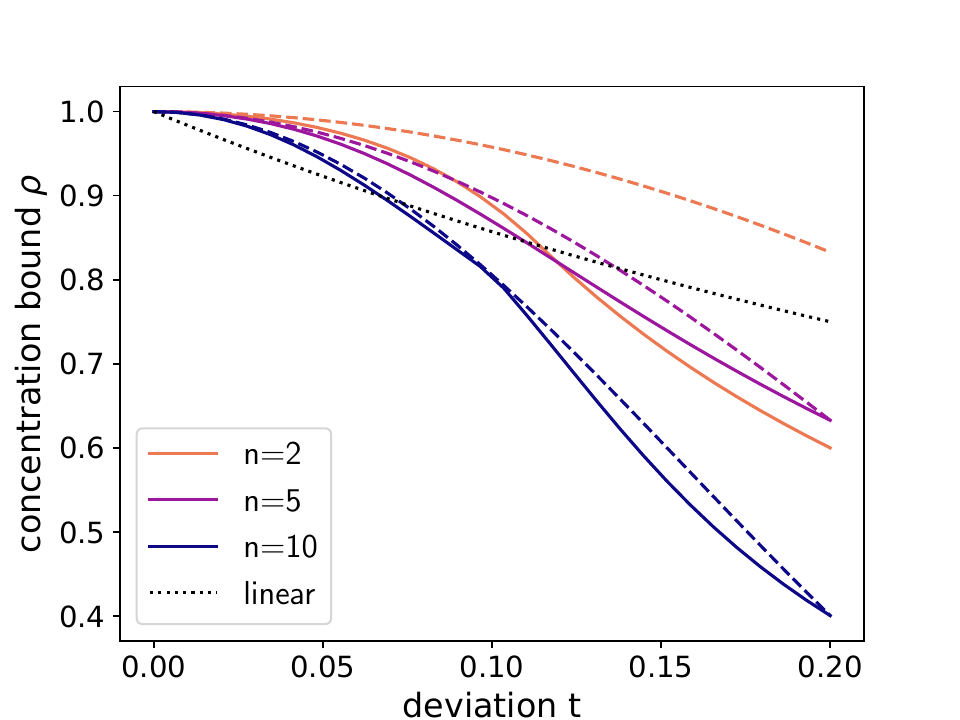}
    \caption{Comparison of the bound obtained in the linear approach~\eqref{eq:linear_upper_optimization problem} to the variational~\eqref{eq:rho_var_conv_iid} and separable~\eqref{eq:exp_multivariate_iid} (in dashed lines) approaches, for $\mu = 0.6$ and several values for $n$.}
    \label{fig:hoeffding_linear_exponential_comparison}
\end{figure}

Surprisingly in Figure~\ref{fig:hoeffding_linear_exponential_comparison}, it happens to improve the variational~\eqref{eq:rho_var_conv_iid} and separable~\eqref{eq:exp_multivariate_iid} approach for small values of the deviation $t$, even for a large number of variables~$n$. However, for larger values of the deviation, the variational and separable approaches significantly outperform the linear approach. In the following, we introduce a family of polynomials encompassing linear functions, that comes at the cost of computationally more expensive approximations.

\subsection{A polynomial upper bound: a SoS approximation for polynomial moment assumptions.}

This section studies an alternate approach to~\eqref{eq:independent_upper_bound_problem} using a family of polynomial upper bounds. Under some assumptions on the degree of these polynomials, tt turns out that this family contains the linear and variational approaches, as feasible points. Its associated Problem~\eqref{eq:research_question_upper_bound} is therefore guaranteed to provide better approximations to the generalized problem of moments. 

We assume $X_1, \ldots, X_n$ to be independent random variables with monomials moments of degree at most $a \in \mathrm{N}^{\star}$, that is $\forall i, \E[g_i(X_i)] = (\E[X_i^k])_{k = 1, \ldots, a} = (\mu_i^{(k)})_{k=1, \ldots, a}$.
Denote $J_d = \{(k_1, \ldots k_n) \in \mathrm{N}^n, \ k_i \in \mathrm{N}, \ k_1 + \ldots + k_n \leqslant d\}$ and $ \forall \kappa \in J_d$, $\bar{x}^\kappa = \prod_{i=1}^n x_i^{k_i}$ monomials of degree $d$. The family of (multivariate) polynomials under consideration is 
\begin{equation}
\label{eq:polynomial_family}
    \mathcal{Q}_{d}^n = \left\{Q(x) = \sum_{\kappa \in J_d} q_\kappa \bar{x}^\kappa, q_\kappa \in \R^{|J_d|}\right\},
\end{equation}
with a certain degree $d \in \mathrm{N}^+$. In particular, linear functions belongs to this set $\mathcal{L} \subset \mathcal{Q}_{d}^n$ along with product functions for $d \geqslant n$. Finally, recall that $F(x) = \mathbf{1}_{x \in S}$, where $S$ has a structure specified in Assumption~\ref{assumption:putinar}.

\begin{assumption}
\label{assumption:putinar}
    Let $S = \{x \in \R^n, h_1(x) \geqslant 0, \ldots, h_m(x) \geqslant 0 \}$ and assume that there exists $s(x)$ a polynomial such that $s(x) = s_0(x) + \sum_{j=1}^m s_j(x) h_j(x)$, with $\{x \in \R^n, h(x) \geqslant 0\}$ a compact set and where $h_i(x)$ are polynomial that admits a sum-of-square decomposition.
\end{assumption}
In particular, any compact polyhedron verifies Assumption~\ref{assumption:putinar}~\citep[Theorem 4.1]{2005Bertsimas_popescu}. Then, Problem~\eqref{eq:research_question_upper_bound} associated with polynomials in $\mathcal{Q}_{d}^n$ takes the form, 
\begin{equation}
\label{eq:ref_polynomial_problem}
\begin{aligned}
    \rho_n^{{\rm polynomial}, d}  =  \inf_{Q \in \mathcal{Q}_{d}^n}& \sup_{\forall i, p_i \in \Px(\X)}  \int_{\X} \cdots \int_{\X} Q(x_1, \ldots, x_n) p_1(x_1) \cdots dp_n(x_n), \\
     {\rm \ such \ that \ } & \forall i, \forall k \in 1, \ldots, a, \int_{\X} x_i^k dp_i(x_i) = \mu_i^{(k)}, {\rm \ and \ } \forall x \in \X^n, \mathbf{1}_{x \in S} \leqslant Q(x). 
\end{aligned}
\end{equation}
By construction, $\rho_n \leqslant \rho_n^{\rm polynomial, a}$.
\begin{proposition}
    \label{prop:reformulation_infinimum_problem}
Let $d \leqslant a$, and let the covariance matrix take the form $\sigma_\kappa = \prod_{i=1}^n \mu_i^{(k_i)}$ for $\kappa \in J_d$. Then Problem~\eqref{eq:ref_polynomial_problem} reformulates as: 
    \begin{equation}
\label{eq:def_polynomial_upper_bound}
\begin{aligned}
    \rho_n^{{\rm polynomial}, d}  =  \inf_{q_\kappa \in \R^{|J_d|}} \sum_{\kappa \in J_d} q_\kappa \sigma_\kappa {\rm \ such \ that \ } \forall x \in \X^n, \mathbf{1}_{x \in S} \leqslant\sum_{\kappa \in J_d} \bar{x}_\kappa q_\kappa.
\end{aligned}
\end{equation}
\end{proposition}
\textbf{Proof. } By definition $Q\in \mathcal{Q}_{d}^n$, $\forall x \in \X^d, Q(x) = \sum_{\kappa \in J_d} q_\kappa \bar{x}^\kappa$. Since $d \leqslant a$ and $\forall i, \forall k \in 1, \ldots, a, \int_{\X} x_i^k dp_i(x_i) = \mu_i^{(k)}$, the objective function $\int_{x \in \X^n}Q(x) dp_1(x_1)\cdots dp_n(x_n)$ for $Q \in \mathcal{Q}_{d}^n$ can be decomposed as the weighted sum of moment $\mu_i^k$. \hfill \Halmos

By construction, Proposition~\ref{prop:reformulation_infinimum_problem} reveals a condition on the degree $d$, so that $\E[Q]$ for $Q \in \mathcal{Q}_d$ formulates as a product combination of moments $(\mu_i^{(k)})_{i, k}$. In addition, Problem ~\eqref{eq:def_polynomial_upper_bound} happens to be exactly the Lagrangian dual of the generalized problem of moments applied to multivariate random variable in $\R^n$~\cite[Equation 2.2]{2005Bertsimas_popescu}. The major difference lies in the structure of the matrix $\forall \kappa \in J_d, \sigma_{\kappa} = \E[\prod_{i=1}^n X_i^{k_i}] = \prod_{i=1}^n E[X_i^{k_i}]$ due to independence.

Problem~\eqref{eq:def_polynomial_upper_bound} suffers from an infinite number of constraints ``$\forall x \in \X^n, \mathbf{1}_{x \in S} \leqslant \sum_{\kappa \in J_d} \bar{x}_{\kappa} q_\kappa$''. Connections between such nonnegative polynomial and sum-of-square decomposition have been hightlighted by several authors~\citep{2002_Lasserre, 2008_Lasserre, 2018deKlerkLaurent}. They were later formulated as SDPs by~\cite{2003_parrilo} and~\cite{2001Lasserre}. In their study of optimal bounds for the generalized problem of moments, \cite{2005Bertsimas_popescu} constructed a sequence of SDPs approximating $\rho_n^{{\rm polynomial}, d}$, that is recalled below.

\begin{theorem}{\cite[Theorem 4.3]{2005Bertsimas_popescu}}
\label{th:sdp_representation}
    Let $S = \{h_j(x) \geqslant 0, j=1, \ldots, m \}$ and $\X = \{\omega_j(x) \geqslant 0, j=1, \ldots, l\}$ verify Assumption~\ref{assumption:putinar}. For every $\epsilon > 0$, there exists a nonnegative integer $r \in \mathrm{N}$ such that $|\rho_n^{{\rm polynomial}, d} - \tilde{\rho}_n^{\rm polynomial, d}(r)| \leqslant \epsilon$, where $\tilde{\rho}_n^{\rm polynomial, d}(r)$ is the value of the following SDP:
    \begin{equation}
\label{eq:sos_representation_polynomial_upper_bound}
    \begin{aligned}
        \tilde{\rho}_n^{\rm polynomial, d}(r) = \inf_{q_\kappa, S \succcurlyeq 0, P \succcurlyeq 0} & \sum_{\kappa \in J_d} q_\kappa \sigma_\kappa, \\
        {\rm \ such \ that \ } & q_\kappa - \delta_{\kappa = 0} = s_\kappa^0 + \sum_{i=1}^m \sum_{\eta, \theta \in J_r, \eta + \theta = \kappa} s_{\eta}^i h_{\theta}^i, \ \forall \kappa \in J_d, \\
       & 0 = \sum_{i=1}^m \sum_{\eta, \theta \in J_r, \eta + \theta = \kappa} s_{\eta}^ih_{\theta}^i, \ \forall \kappa \in J_r \backslash J_d, \\
       & q_\kappa  = p_\kappa^0 + \sum_{i=1}^l \sum_{\eta, \theta \in J_r, \eta + \theta = \kappa} p_{\eta}^i \omega_{\theta}^i, \ \forall \kappa \in J_d, \\
       & 0  = p_\kappa^0 + \sum_{i=1}^l \sum_{\eta, \theta \in J_d, \eta + \theta = \kappa} p_{\eta}^i \omega_{\theta}^i, \ \forall \kappa \in J_r \backslash J_d, \\
       & s_\kappa^i = \sum_{\eta, \theta \in J_r, \eta + \theta = \kappa} s_{\eta, \theta}^i, \ S^i = [s_{\eta, \theta}]_{\eta, \theta \in J_r} \succcurlyeq 0,\ \forall \kappa \in J_r, \ i=1, \ldots, m, \\
       & p_\kappa^i = \sum_{\eta, \theta \in J_r, \eta + \theta = \kappa} p_{\eta, \theta}^i, \ P^i = [p_{\eta, \theta}]_{\eta, \theta \in J_r} \succcurlyeq 0, \  \forall \kappa \in J_r, \ i=1, \ldots, l.
    \end{aligned}
\end{equation}
\end{theorem}

Theorem~\ref{th:sdp_representation} provides a sequence of SDPs approximation Problem~\eqref{eq:def_polynomial_upper_bound}, but does not specify at which degree $r$ a certain precision level $\epsilon$ is attained. By construction, it is clear that their degree $r$ must be larger than the degree $d$ of the polynomials under consideration. \cite{2005Bertsimas_popescu} highlighted a hierarchy between these SDP approximations, which is explicited in Corollary~\ref{corollary:hierarchie_sdp}.
\begin{coro}{\cite[Corollary 4.4]{2005Bertsimas_popescu}}
\label{corollary:hierarchie_sdp}
Let $\rho_n^{{\rm polynomial}, d}$ be defined as in~\eqref{eq:def_polynomial_upper_bound} and $\tilde{\rho}_n^{ \rm polynomial, d} (r)$ as in~\eqref{eq:sos_representation_polynomial_upper_bound}. Then, it follows that:
    \begin{equation}
        \rho_n \leqslant \rho_n^{{\rm polynomial}, d} \leqslant \tilde{\rho}_n^{\rm polynomial, d} (r) \leqslant \cdots \leqslant \tilde{\rho}_n^{\rm polynomial, d} (1).
    \end{equation}
\end{coro}
From a numerical perspective, the size of SDPs in~\eqref{eq:sos_representation_polynomial_upper_bound} grows exponentially with degrees $r$ and $d$, limiting both the precision of the approximation problem~\eqref{eq:def_polynomial_upper_bound}. In the next section however, we apply these approximations to two basic concentration inequalities and improve the existing bounds for a small number of variables $n$ and a small degree $r$.

\subsection{Applications to Bernstein's and Bennett's inequalities}

The variational and separable approaches failed to provide tractable optimization problem for second-order conditions, that appear for example in Bennett's and Bernstein's inequalities. For two random variables, we compute refined bounds using the polynomial approach, and more precisely the SDP aproximations defined in~\eqref{eq:sos_representation_polynomial_upper_bound}.

\textbf{Bernstein's inequality.} Bernstein's inequality controls the deviation of the sum of independent random variables to their mean, given an appropriate control of moments. There exists different versions for Bernstein's inequality, and we consider a convenient version requiring finite second-order moments, as provided for example in~\cite[Corollary 2.11]{2013Boucheron}.

\begin{corollary}
\label{coro:bernstein}
    Let $X_1, \ldots, X_n$ be independent random variables in $[-c, c]$ with means $\E[X_i] = \mu_i$, and variance $v = \sum_{i=1}^n \E[X_i^2]$. Then, for all $t > 0$, 
    \begin{equation}
    \label{eq:bernstein_coro}
        \Prob\left(\sum_{i=1}^n X_i \geqslant nt + \sum_{i=1}^n \mu_i\right) \leqslant \exp\left(- \frac{n^2 \frac{t^2}{2}}{v + c \frac{nt}{3}} \right).
    \end{equation}
\end{corollary}
Let $X_1, X_2$ be two i.i.d. random variables  taking their values almost surely in $[-1, 1]$ with $\E[X_1] = \E[X_2] = \mu^{(1)}$, and $\E[X_1^2] = \E[X_2^2] = \mu^{(2)}$. We consider the polynomial family $\mathcal{Q}_2 = \{ x \in \R^2 \mapsto x_{(2)}^\top Q x_{(2)}, Q \in \R^{(6 \times 6)}\}$ where $x_{(2)} = (1, x_1, x_2, x_1x_2, x_1^2, x_2^2)$. Then, Problem~\eqref{eq:def_polynomial_upper_bound} takes the form: 
\begin{equation}
\label{eq:polynomial_n2_bernstein}
\begin{aligned}
    \forall t \in \geqslant 0, \ \rho_2^{\rm polynomial, 2}(t) = \inf_{Q \in \R^{(6 \times 6)} }  \  {\rm Tr}(Q \Sigma) \ {\rm such \ that \ }  & \forall x \in \X^2, \ 1 \leqslant x_{(2)}^\top Q x_{(2)}, \\
& \forall x \in [-1, 1]^2, 0 \leqslant x_{(2)}^\top Q x_{(2)},
\end{aligned}
\end{equation}
where $\bar{\X}_2 = \{(x_1, x_2) \in [-1, 1], x_1 + x_2 \geqslant 2(\mu^{(1)} + t)\}$, $\Sigma = \sigma \sigma^\top$ and $\sigma = (1, \mu^{(1)}, \mu^{(1)}, (\mu^{(1)})^2, \mu^{(2)}, \mu^{(2)})$.  It results from Theorem~\ref{th:sdp_representation} that Problem~\eqref{eq:polynomial_n2_bernstein} can be approximated by a hierarchy of sum-of-square~\eqref{eq:sos_representation_polynomial_upper_bound}. 

\begin{figure}[h!]
     \centering
     \begin{subfigure}[t]{0.45\textwidth}
         \centering
         \includegraphics[height=120pt]{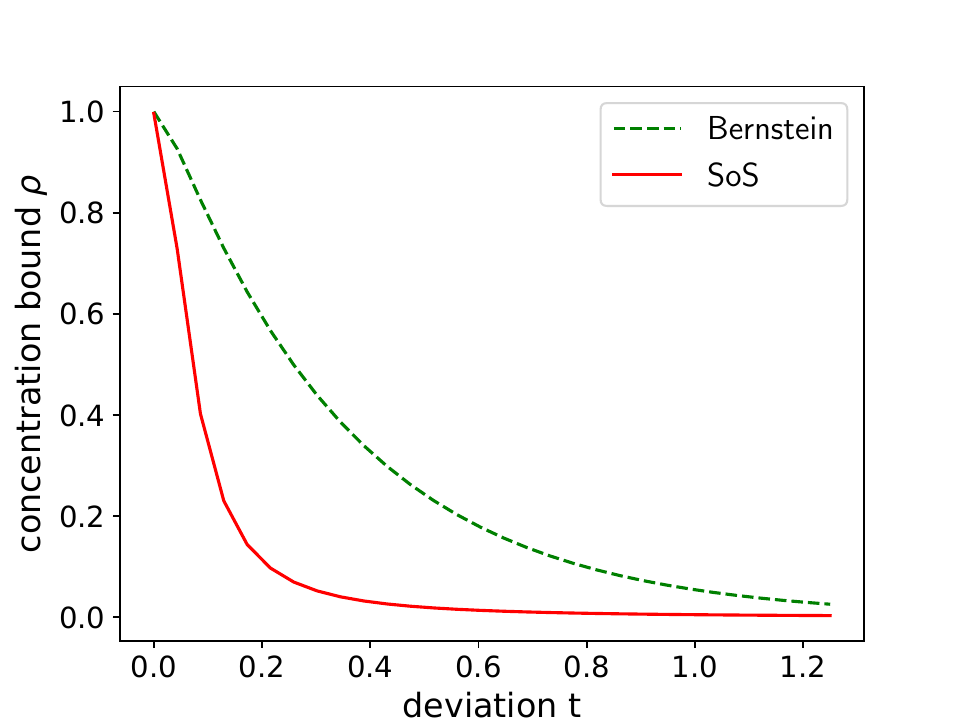}
         \caption{$\mu_1 = -0.3$, $\mu_2 = 0.1$}
     \end{subfigure}
     \hfill
     \begin{subfigure}[t]{0.45\textwidth}
         \centering
         \includegraphics[height=120pt]{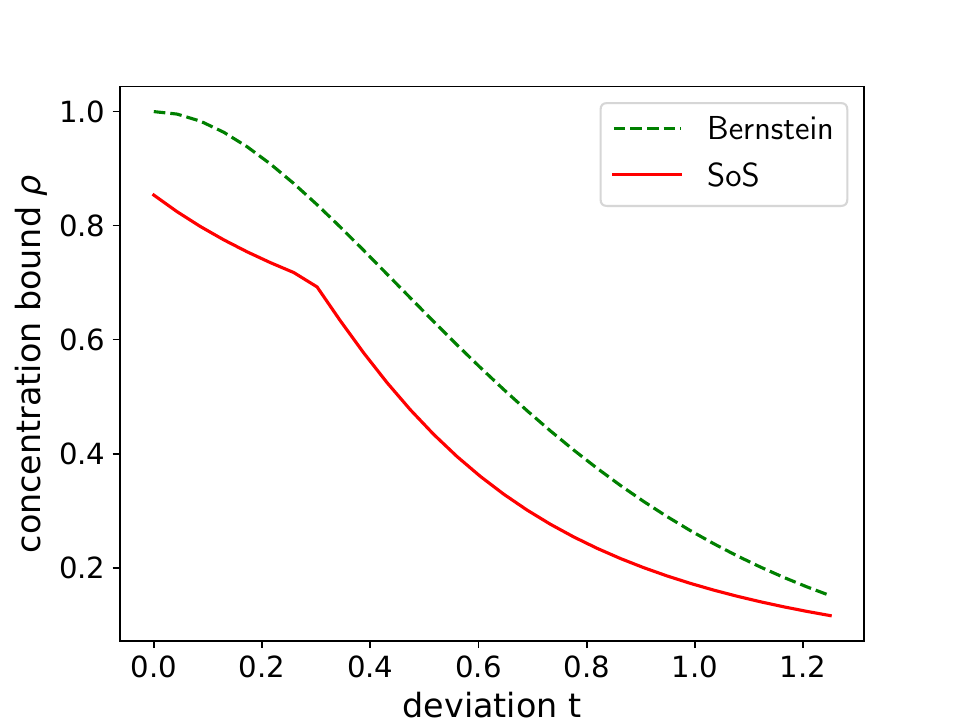}
         \caption{$\mu_1 = -0.3$, $\mu_2 = 0.5$}
     \end{subfigure}
     \caption{\normalsize{Comparison of Bernstein's inequality~\eqref{eq:bernstein_coro} to an SDP approximation~\eqref{eq:sos_representation_polynomial_upper_bound} to the polynomial approach~\eqref{eq:polynomial_n2_bernstein}, as a function of $t$, for $n=2$, $d=2$ and $r=2$.}}
\label{fig:bernstein_function_t}
\end{figure}
Figure~\ref{fig:bernstein_function_t} shows that the SDP approximation of the polynomial approach~\eqref{eq:polynomial_n2_bernstein}
improves upon Bernstein's inequality~\eqref{eq:bernstein_coro} for a small degree in the SoS hierarchy (here $r=2$). This improvement is not straightforward, since proofs for Bernstein's inequality often require a specific control of moments $\E[X_i^k]$ for all $k~\in~\mathrm{N}$ (see, e.g.,\cite[proof for Theorem 2.9]{2013Boucheron}). Therefore, considering higher-degree $d$ in the polynomial approximations would probably produce a lower value for $\rho_2^{\rm polynomial, d}$. However, increasing $d$ requires to increase the degree $r \in \mathrm{R}$ of the sum-of-square approximations~\eqref{eq:sos_representation_polynomial_upper_bound} and thereby, leads to very large SDPs that cannot be handled by our solvers. 

\textbf{Bennett's inequality.} Bennett's inequality is similar with Bernstein's inequality, but applies to upper bounded random variables, as recalled in Theorem~\ref{th:bennet}.

\begin{theorem}{\cite{1968Bennet}}
\label{th:bennet}
    Let $X_1, \cdots, X_n$ be i.i.d. random variables such that $X_i \leqslant a$ almost surely and $\sigma^2 = \sum_{i=1}^n \E[(X_i - \E X_i)^2]$, then for any $t\geqslant 0$,
    \begin{equation}
    \label{eq:bennet_th}
        \Prob\left(\sum_{i=1}^n \{ X_i - \E[X_i]\} \geqslant nt\right) \leqslant {\rm exp}\left( - \frac{\sigma^2}{a^2}h(\frac{atn}{\sigma^2}) \right),
    \end{equation}
    where $h(t) = (1 + t) \log{(1 + t)} - t$. It implies that $\Prob(\sum_{i=1}^n \{ X_i - \E[X_i]\} \geqslant nt) \leqslant \exp{\left(-\frac{t^2n^2}{2(\sigma^2 + atn/3)}\right)}$.
\end{theorem}

Let $X_1, X_2$ be two i.i.d. random variables  taking their values almost surely in $\mathcal{Z} = [-\infty, 1]$ with $\E[X_1] = \E[X_2] = \mu^{(1)}$, and $\E[X_1^2] = \E[X_2^2] = \mu^{(2)}$. Then, Problem~\eqref{eq:def_polynomial_upper_bound} takes the form: 
\begin{equation}
\label{eq:polynomial_n2_bennet}
\begin{aligned}
    \rho_2^{\rm polynomial, 2}(t) = \inf_{Q \in \R^{6 \times 6}}  \  {\rm Tr}(Q \Sigma) \ {\rm such \ that \ }  & \forall x \in \bar{\mathcal{Z}}_2, \ 1 \leqslant x_{(2)}^\top Q x_{(2)}, \\
& \forall x \in [-\infty, 1]^2, 0 \leqslant x_{(2)}^\top Q x_{(2)},
\end{aligned}
\end{equation}
where $\bar{\mathcal{Z}}_2 = \{(x_1, x_2) \in [-\infty, 1], x_1 + x_2 \geqslant 2(\mu^{(1)} + t)\}$, $\Sigma = \sigma \sigma^\top$ and $\sigma = (1, \mu^{(1)}, \mu^{(1)}, (\mu^{(1)})^2, \mu^{(2)}, \mu^{(2)})$.

\begin{figure}[h]
     \centering
     \begin{subfigure}[t]{0.45\textwidth}
         \centering
         \includegraphics[height=120pt]{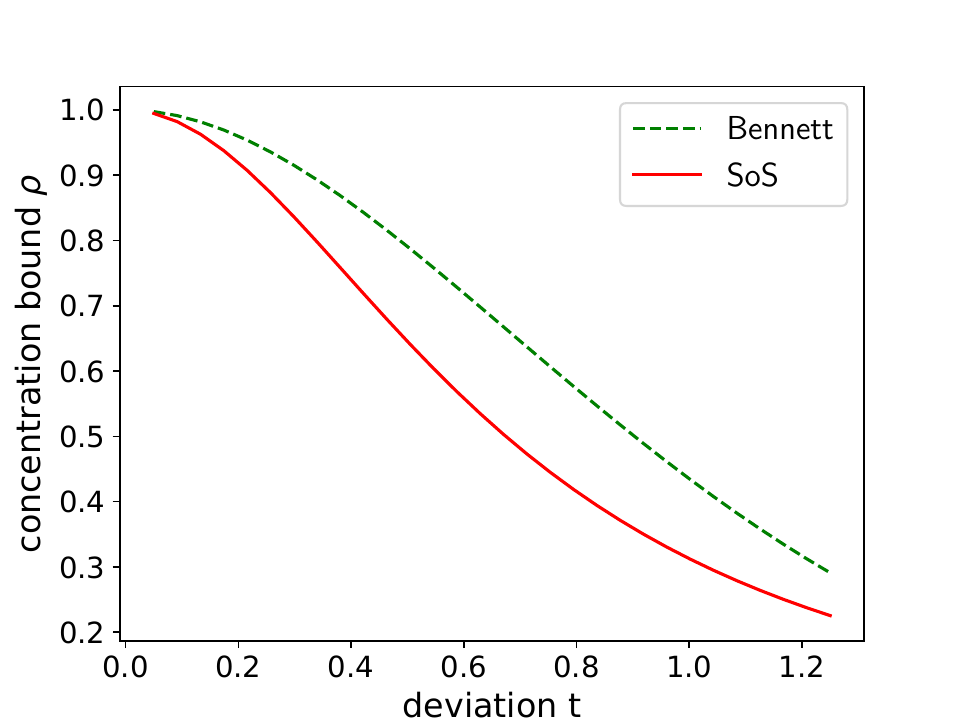}
         \caption{$\mu_1 = -0.3$, $\mu_2 = 1$}
     \end{subfigure}
     \hfill
     \begin{subfigure}[t]{0.45\textwidth}
         \centering
         \includegraphics[height=120pt]{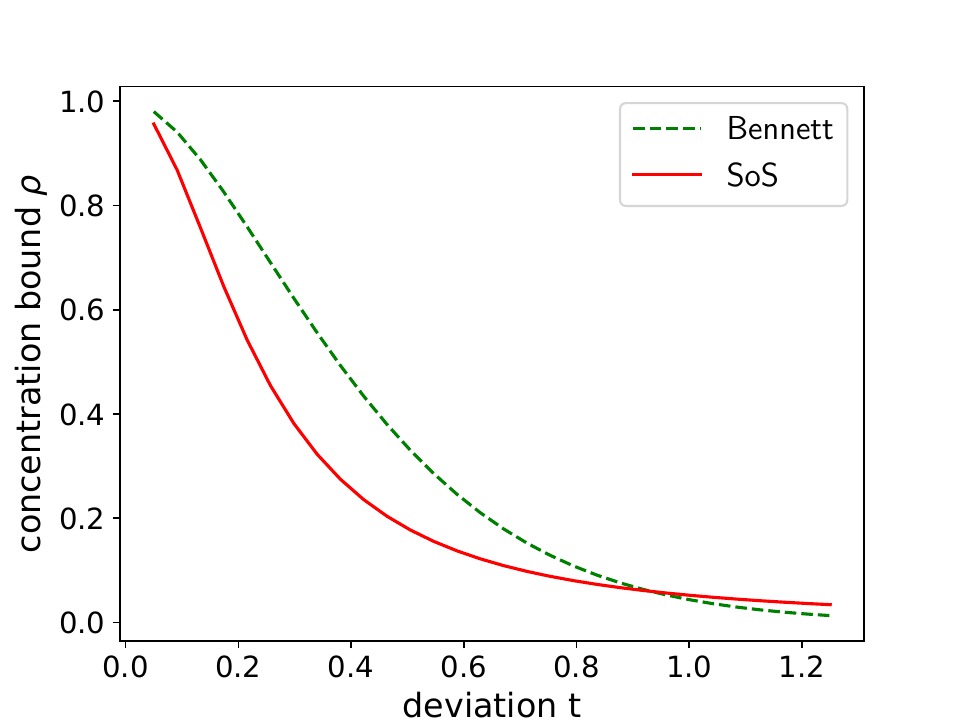}
         \caption{$\mu_1 = -0.3$, $\mu_2 = 0.2$}
     \end{subfigure}
     \caption{\normalsize{Comparison of Bennett's inequality~\eqref{eq:bennet_th} to an SDP approximation~\eqref{eq:sos_representation_polynomial_upper_bound} to the polynomial approach~\eqref{eq:polynomial_n2_bennet}, as a function of $t$ for $n=2$, $d=2$ and $r=2$.}}
\label{fig:bennet_function_t}
\end{figure}

Figure~\ref{fig:bennet_function_t} reveals the limit of the polynomial approach, which fails to improve Bennett's inequality for all values of $t \geqslant 0$. The proof relies indeed on additional arguments, such as Jensen's inequality, or Taylor approximations such as $u \mapsto \log(1 + u) \leqslant u$, which are not exploited in this approach. 

To conclude, the polynomial approach~\eqref{eq:def_polynomial_upper_bound} effectively refines Bernstein and Bennett's inequalities for some range of values $t \geqslant 0$. However it is hindered by the necessity of large SDP approximations: handling higher-order polynomials and a large number of random variables increases the size of the SDPs under consideration, making them difficult to solve. A review of the complexity of semidefinite optimization and interior point methods can be found in~\cite{Wolkowicz2000, 1996Vandenberghe}.

\subsection{Higher-order polynomial approximations : a variational reformulation.}

The variational, separable and linear approaches allowed refining Hoeffding's inequality. In the polynomial approach developed above, Proposition~\ref{prop:reformulation_infinimum_problem} entails that the degree of the polynomials under consideration is controlled by the  highest-order moment. In the context of Hoeffding's inequality, where the first-order moment is finite, it implies that the best polynomial representation is actually linear ($d \leqslant 1$).

In the following, we introduce a numerical procedure to incorporate higher-order polynomials for studying Hoeffding's inequality. We then lay the foundation for a feature-based approach that generalizes the polynomial approach to broader families of upper functions. Finally, we derive a SDP relaxation for the case of two independent random variables.

\subsubsection{Studying Hoeffding's inequality with second-order polynomials.}

Hoeffding's inequality requires finite first-order moments, together with almost surely bounded variables. In what follows, we show how it affects second-order moments. Based on this, we integrate second-order moment conditions into a polynomial approach for studying Hoeffding's inequality.

Let $X_1, \ldots, X_n$ be independent random variables taking their values in $[0, 1]$ almost surely with $\E[X_i] = \mu_i^{(1)}$ for all $i=1, \ldots, n$. Lemma~\eqref{proof_lem:var_to_mean_rel} ensures the existence and a control on the second-order moment (see proof in Appendix~\ref{proof_lem:var_to_mean_rel}).

\begin{lemma}
\label{lem:var_to_mean_relationship}
    Let $X$ be a random variable almost surely on $[0, 1]$ with mean $\E[X] = \mu^{(1)}$.
    Then $X$ admits a finite second-order moment, denoted $\E[X^2] = \mu^{(2)}$. In addition, it holds that $(\mu^{(1)})^2 \leqslant \mu^{(2)} \leqslant \mu^{(1)}.$
\end{lemma}

\begin{figure}[h]
     \centering
     \begin{subfigure}[t]{0.32\textwidth}
         \centering
         \includegraphics[height=120pt]{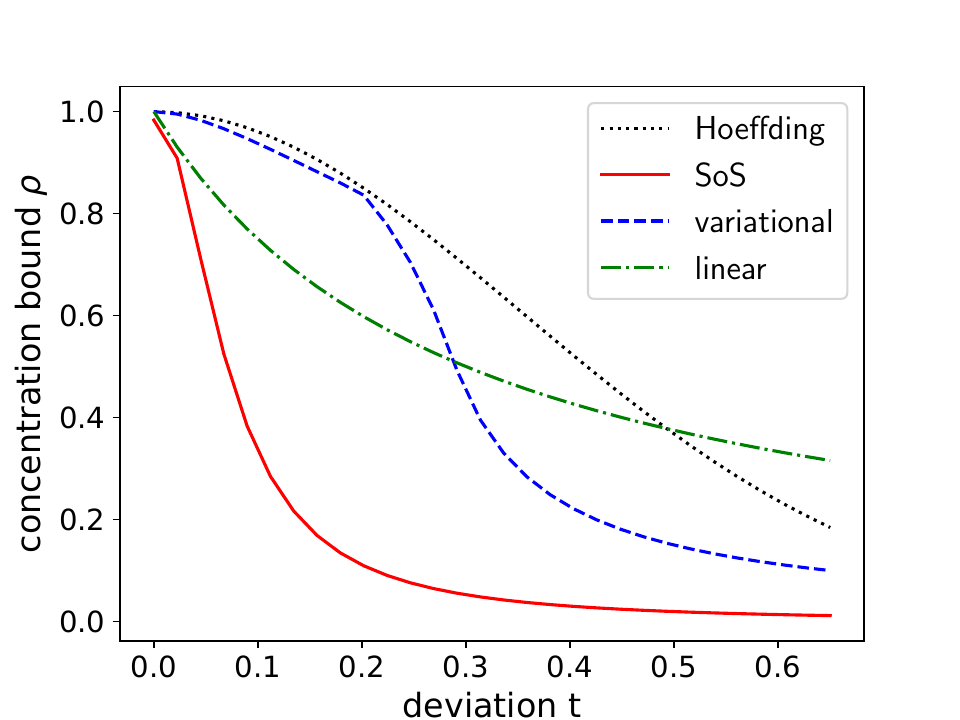}
         \caption{$\mu_1 = 0.3$ and $\mu_2 = 0.1$}
     \end{subfigure}
     \hfill
     \begin{subfigure}[t]{0.32\textwidth}
         \centering
         \includegraphics[height=120pt]{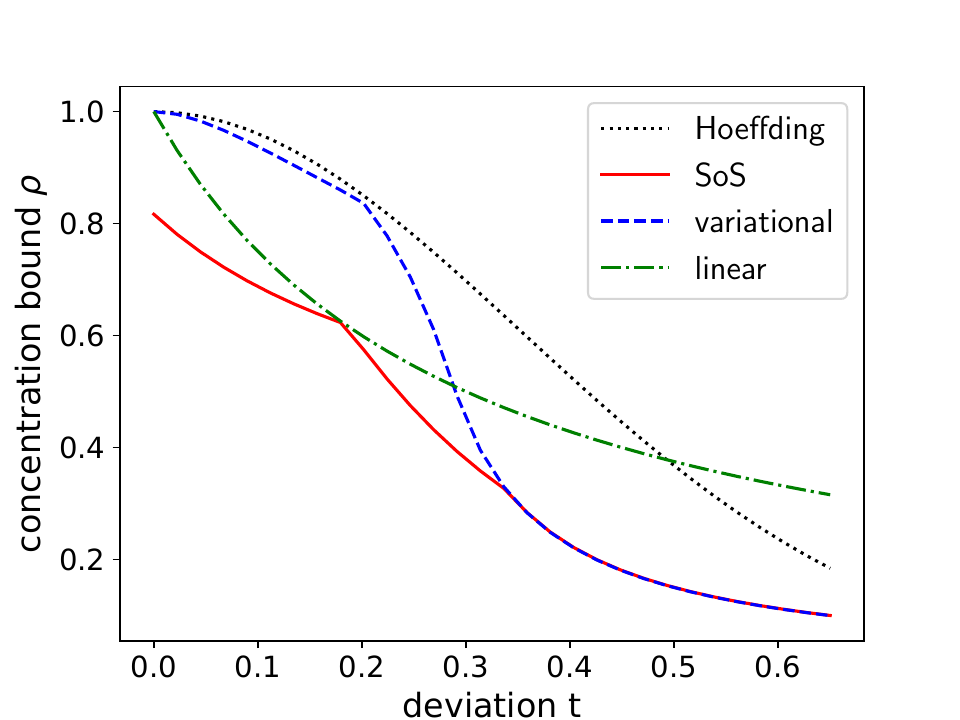}
         \caption{$\mu_1 = 0.3$ and $\mu_2 = 0.2$}
     \end{subfigure}
     \hfill
    \begin{subfigure}[t]{0.32\textwidth}
    \centering
    \includegraphics[height=120pt]{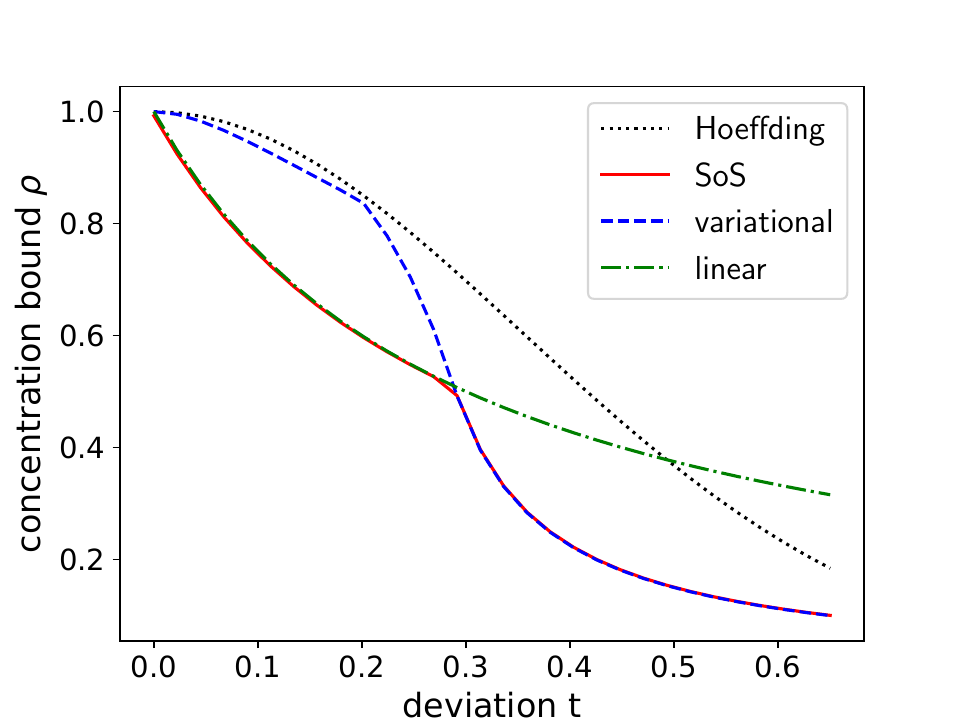}
    \caption{Optimized version for $\mu_1 = 0.3$.}
\end{subfigure}
     \caption{\normalsize{Comparison of the bound in polynomial approach~\eqref{eq:def_polynomial_upper_bound} (with an SDP approximation of degree $r=3$), to the variational~\eqref{eq:rho_var_conv_iid} and linear~\eqref{eq:linear_upper_optimization problem} approaches, and to Hoeffding's inequality~\eqref{eq:Hoeffding_value}, as a function of the deviation parameter $t$ for $n=2$.}}
\label{fig:sos_several_t}
\end{figure}

Lemma~\ref{lem:var_to_mean_relationship} provides a control of the second-order moment by the first-order moment for bounded random variables. From that, we define an optimization problem based on the polynomial approach~\eqref{eq:def_polynomial_upper_bound}:
\begin{equation}
\label{eq:numerical_procedure_optimization_mu2}
\begin{aligned}
    \tilde{\rho}_n^{\rm polynomial, 2}  &=  \inf_{\mu^{(2) \in [(\mu^{(1)})^2, \mu^{(1)}]}} \ \inf_{Q \in \mathcal{Q}_{2}^n} \sup_{p_1, \ldots, p_n \in \Px(\X)} & \int_{\X} \ldots \int_{\X}  Q(x_1, \ldots, x_n) dp_1(x_1) \cdots dp_n(x_n), \\
   &  & {\rm \ such \ that  \ } \forall i, \int_{\X} x_i dp_i(x_i) = \mu_i^{(1)}, \\
   &  & \forall i, \int_{\X} x_i^2 dp_i(x_i) = \mu_i^{(2)}, \\
   &  & \forall x \in \X^n, \mathbf{1}_{ x \in S} \leqslant Q(x), \\
  \tilde{\rho}_n^{\rm polynomial, 2}  &= \inf_{\mu^{(2) \in [(\mu^{(1)})^2, \mu^{(1)}]}} \rho_n^{\rm polynomial, 2}(\mu^{(2)}).
\end{aligned}
\end{equation}
Given a value of $\mu^{(2)} \in[(\mu^{(1)})^2, \mu^{(1)}]$, the inner optimization problem defining $\rho_n^{\rm polynomial, 2}$ can be approached with a hierarchy of sum-of-square defining $\rho_n^{\rm polynomial, 2}(r, \mu^{(2)})$, as stated in Proposition~\ref{prop:reformulation_infinimum_problem} and Theorem~\ref{th:sdp_representation}. In Figure~\ref{fig:sos_several_t}, we use a gridsearch procedure on $\mu^{(2)}$ for approximating $\rho_n^{\rm polynomial, 2}(\mu^{(2)})$. It turns out that optimizing $\rho_n^{\rm polynomial, 2}(\mu^{(2)})$ with respect to $\mu^{(2)}$ aligns exactly with the minimum bound achieved by the linear~\eqref{eq:linear_upper_optimization problem} and variational~\eqref{eq:rho_var_conv_iid} approaches. 

In conclusion, linear functions and product-functions appear to be the minimal polynomial families suitable for studying Hoeffding's inequality. Although the polynomial approach~\eqref{eq:numerical_procedure_optimization_mu2} incurs high computational costs, the linear approach offers a closed-form solution and the variational approach provides tractable reformulations for i.i.d. random variables and structured independent variables.

\subsubsection{Generalization to tighter bounds for $F$: a variational feature-based formulation}
We have previously explored the family of product-functions, linear functions as well as polynomials. These approaches were limited by the power of representation allowed by moments (see Proposition~\ref{prop:reformulation_infinimum_problem}). In the context of Hoeffding's inequality, we addressed higher-order polynomials using finite lower-order moments through a numerical procedure. We now extend this approach by laying the foundations for a feature-based approach.

 Let us start by recalling the variational formulation~\eqref{eq:research_question_upper_bound}:
\begin{equation*}
    \inf_{H \in \mathcal{H}} \sup_{\forall i, p_i \in \Px_{\mu_i}(\X)} \int_{\X^n}H(x)dp_1(x_1) \cdots dp_n(x_n) {\rm \ such \ that \ } \forall x \in \X^n, F(x) \leqslant H(x),
\end{equation*}
 where $\mathcal{H}$ is a family of upper bounds. We introduce a feature vector $\phi:\X \mapsto \R^l, l\in \mathrm{N}$ such that functions have the following representation:
\begin{equation*}
    F(x) = \Biggl\langle \bar{F}, \bigotimes_{i=1}^n \phi(x_i)\biggr \rangle = \sum_{i_1, \ldots, i_n = 1}^l \bar{F}_{i_1, \ldots, i_n} \prod_{k=1}^n\phi(x_k)_{i_k}.
\end{equation*}
In the polynomial (resp.~variational) approach for examples, features $\phi$ were monomials (resp.~linear functions). Let us consider a family of functions $H$ decomposing with respect to these features:
\begin{equation*}
\begin{aligned}
    \sup_{\forall i, p_i \in \Px_{\mu_i}(\X)} \int_{\X^n}H(x)dp_1(x_1) \cdots dp_n(x_n) &= \sup_{\forall i, p_i \in \Px_{\mu_i}(\X)} \Biggl\langle \bar{H}, \bigotimes_{i=1}^n\int_\X\phi(x_i) dp_i(x_i)\biggl \rangle = \sup_{\forall i, \sigma_i \in \mathcal{K}(\mu_i)} \Biggl\langle \bar{H}, \bigotimes_{i=1}^n \sigma_i \biggr \rangle,
\end{aligned}
\end{equation*}
where $\mathcal{K}(\mu_i)$ is the set of achievable moments $\E_{p_i \in \mathcal{P}_{\mu_i}(\X)}[\phi(X_i)] = \int_\X \phi(x_i)dp_i(x_i)$ such that $p_i \in \Px_{\mu_i}(\X)$. In the context of Hoeffding's inequality, the set of achievable moments~$\mathcal{K}(\mu_i)$ corresponds to the intuitive idea that there is no assumption on second-order moments, but that they are related to lower-order moments. Then, the overall optimization problem takes the form: 
\begin{align}
\label{eq:feature_based_reformulation}
    \inf_{H \in \mathcal{H}} \sup_{\forall i, \sigma_i \in \mathcal{K}(\mu_i)} \Biggl\langle \bar{H}, \bigotimes_{i=1}^n \sigma_i \biggr\rangle {\rm \ such \ that \ } \forall x \in \X^n, \Biggl\langle \bar{F} - \bar{H}, \bigotimes_{i=1}^n \phi(x_i) \biggr\rangle \leqslant 0.
\end{align}
Compared to the polynomial approach developed above~\eqref{eq:def_polynomial_upper_bound}, Problem~\eqref{eq:feature_based_reformulation} provides a more generic formulation for any family of features. The family of features must satisfy two key components : a (tighter) relaxation of the constraint $\forall x \in \X^n, \langle \bar{F} - \bar{H}, \bigotimes_{i=1}^n \phi(x_i) \rangle \leqslant 0$, which can be managed using sum-of-square for polynomial features, and a relaxation for $\sup_{\forall i, \sigma_i \in \mathcal{K}(\mu_i)} \langle \bar{H}, \bigotimes_{i=1}^n \sigma_i \rangle$. In the following proposition, we analyze a simple relaxation of $\sup_{\forall i, \sigma_i \in \mathcal{K}(\mu_i)} \langle \bar{H}, \bigotimes_{i=1}^n \sigma_i \rangle$ for two i.i.d. random variables in the context of Hoeffding's inequality (i.e., $n=2$). For this case, the tensor formulation simplifies into matrices and admits an SDP relaxation.

\begin{proposition}
\label{prop:two_iid_relaxation}
    Let $X_1, X_2$ be i.i.d. random variables taking their values almost surely in $[0, 1]$ with finite mean $\mu^{(1)}$ and let $\phi(x) = (1, x, x^2)$ be the feature vector. Then, 
    \begin{align*}
        \sup_{\sigma \in \mathcal{K}(\mu^{(1)})} {\rm Tr}\left( \tilde{H} \sigma \sigma^\top\right) \leqslant \sup_{M \succcurlyeq 0} \   {\rm Tr}\left( \tilde{H} M \right) {\rm \ such \ that \ } & 0 \leqslant {\rm Tr}(ME_{1, 1}) \leqslant  {\rm Tr}(ME_{0, 2}) \leqslant  {\rm Tr}(ME_{0, 1}) \leqslant 1, \\
        &  {\rm Tr}(ME_{0, 0}) = 1.
    \end{align*}
\end{proposition}
\textbf{Proof.} Let $\tilde{H} \in \R^{3 \times 3}$ be a symmetric matrix representation of vector $\bar{H} \in \R^6$, such that $\sup_{\sigma \in \mathcal{K}(\mu^{(1)})} {\rm Tr}( \tilde{H} \sigma \sigma^\top) = \sup_{M \succcurlyeq 0} {\rm Tr}( \tilde{H} M) {\rm \ such \ that \ } M \in {\rm hull}\{\sigma \sigma^\top , \sigma \in \mathcal{K}(\mu)\}$. By definition, $\sigma = (1, \mu^{(1)}, \mu^{(2)})$ holds for the first and second-order moments. We relax this problem by optimizing over $M \succcurlyeq 0$ and incorporating additional constraints on $M$ to ensure it accurately incorporates relationships between the first and second-order moments (Lemma~\ref{lem:var_to_mean_relationship}, namely $(\mu^{(1)})^2 \leqslant \mu^{(2)} \leqslant \mu^{(1)}$ and $\int_0^1 dp(x) = 1$.). \hfill \Halmos.

Proposition~\ref{prop:two_iid_relaxation} provides an SDP relaxation for approximating Hoeffding's inequality for two random variables and given a quadratic features. Despite the simple relaxation, solving this problem requires to solve efficiently a saddle-point problem and to quantify how far such a relaxation is from the generalized problem of moments~\eqref{eq:independent_upper_bound_problem}, which are left to future work.


\section{Conclusion and future works.} \label{sec:con}

\paragraph{Conclusion.}In this work, we introduced two families of upper bounds for approximating the generalized problem of moments for independent random variables. First, we studied a separable approach, leveraging specific upper functions adapted to finite first-order moments. This approach is complemented by a convex variational method, optimizing over the entire family of product-functions. When studying Hoeffding's inequality, these formulations are particularly effective for both i.i.d. random variables and cases where variables are divided into groups with different means, facilitating the reconstruction of associated extremal distributions. Due to the computational limitations of the product-functions based approaches, we broadened our scope by introducing a polynomial family of upper bounds. Here, carefully selected polynomials are employed, resulting in non-negative polynomials that can be approximated using sum-of-square techniques, at a higher computational cost. This framework enables exploration of concentration properties concerning Bennett's and Bernstein's inequalities, although it does not universally improve theoretical bounds. We finally extended the polynomial approach into a feature-based approach, using polynomials independently of the order of moment assumptions. While not focused on computational efficiency, this method offers a more flexible and comprehensive way to study concentration inequalities. In summary, our methodologies each address different complexities inherent in the problem of moments and independence.

\paragraph{Future works.}
Throughout this work, we have highlighted several limitations related to these approaches. The separable approach could benefit from exploring and constructing new product-functions that lead to closed-form solutions, going beyond moment-generating functions. Meanwhile, the polynomial approach, when approximated by SDPs of increasing sizes, converges slowly even for a reasonable number of variables and low-degree polynomials. Studying the connection between the family of polynomials, the linear functions and product-functions would simplify the underlying optimization problems. In addition, numerical experiments have shown the limitations of the polynomial approach, which does not account for some analytical arguments, such as convexity or inequalities derived from Taylor expansion. Introducing key components of these analyses would probably help improve (or match) known bounds. Furthermore, the feature-based approach could still benefit from exploring efficient relaxations to improve known bounds.

Finally, the generalized problem of moments requires few assumptions on the distribution under consideration, such as bounded moments or random variables lying in bounded sets (e.g., intervals). Other paths to improvements of concentration bounds without independence were explored, by exploiting additional structural properties of the distributions. For instance, \cite{2005_Popescu} extended the generalized problem of moments to convex classes of distributions, such as unimodal or bimodal distributions. Given unimodal distributions, \cite{2015_VanParys_Goulart_Kuhn} reconstructed non-discrete extremal distributions, which are thereby more representative of the effective behaviors of random variables. Extending their work to the case of independent variables would probably help refine inequalities that can be adapted to different problem structures. In particular, introducing subgaussian assumptions, which appear in many probability proofs in machine learning, would be of interest, thus improving probabilistic bounds for such problems.

\begin{appendices}

 \section{Exact formulation of generalized problem of moments for independent random variables using optimal transport: }
 \label{app:optimal_transport}
This proof is an alternative to the convexification of the variational approach from Proposition~\ref{prop:var_approach_convexified}. Using optimal transport, we compute an exact formulation for the generalized problem of moments~\eqref{eq:independent_upper_bound_problem}, and recover the convex reformulation for the variational approach as result of weak duality.

Let us introduce for $x \in \X^n$, $G(x) = \log(F(x_1, ..., x_n))$. Then, the exact formulation takes the form for the KL divergence $D(q\bigg|\bigg|p) = \int_\X \log\left( \frac{p(x)}{q(x)}\right)dq(x)$:
\begin{align*}
    \log(\rho_n) &= \sup_{\forall i, p_i \in \Px_{\mu_i}(\X)} \log\left(\int_{\X^n} e^{G(x)} dp_1(x_1) \cdots dp_n(x_n) \right), \\
    &= \sup_{q \in \Px(\X^n)} \sup_{\forall i, p_i \in \Px_{\mu_i}(\X)} \int_{\X^n} G(x) dq(x) - D\left(q\bigg|\bigg|\prod_{i=1}^n p_i\right),
\end{align*}
By Donsker-Varadhan's inequality. By Pythagorean theorem for the KL divergence and mutual information, we have:
\begin{align*}
    \log(\rho_n) &= \sup_{q \in \Px(\X^n)} \sup_{\forall i, p_i \i \Px_{\mu_i}(\X)} \int_{\X^n} G(x) dq(x) - D\left(q\bigg|\bigg|\prod_{i=1}^n q_i\right) - \sum_{i=1}^n D(p_i \bigg|\bigg|q_i).
\end{align*}
The variational relaxation corresponds in fact to using the fact that $D(q\bigg|\bigg|\prod_{i=1}^n q_i) \geqslant 0$. In addition, we recover the convex formulation from the Lagrangian relaxation of:
\begin{align*}
    \log(\rho_n) = \sup_{q \in \Px(\X^n)} &\inf_{\forall i, \alpha_i \in \R, \beta_i \in \R^m}  \int_{\X^n} G(x) dq(x) - D\left(q\bigg|\bigg|\prod_{i=1}^n q_i\right) + \sum_{i=1}^n \{\alpha_i + \beta_i^\top \mu_i - 1\} \\
 & + \sup_{x \in \X^n } \{G(x) - \sum_{i=1}^n \log(\alpha_i + \beta_i^\top g(x_i))\}, \\
     \leqslant \inf_{\forall i, \alpha_i \in \R, \beta_i \in \R^m} &  \sum_{i=1}^n \{\alpha_i + \beta_i^\top \mu_i - 1\} + \sup_{x \in \X^n} \{G(x) - \sum_{i=1}^n \log(\alpha_i + \beta_i^\top g(x_i))\} \\
     & +\sup_{q \in \Px(\X^n)} \int_{\X^n} G(x) dq(x) - D\left(q\bigg|\bigg|\prod_{i=1}^n q_i\right),
\end{align*}
by weak duality. Ignoring the term, $D\left(q \bigg|\bigg| \prod_{i=1}^n q_i\right)$, the maximum of $\sup_{\X^n} G(x) dq(x) = \max_{x \in \X^n} G(x)$ with respect to $q \in \Px(\X^n)$ is attained at a Dirac, at which $D\left(q \bigg|\bigg| \prod_{i=1}^n q_i\right) = 0$. Thus, the fact that $\rho_n \leqslant  \rho_n^{{\rm var}}$ holds by weak duality.

\section{Proofs for Hoeffding's inequality}
\subsection{Separable approach for univariate distributions}

\subsubsection{Proof for Proposition~\ref{prop:univariate_exponential}}
\label{app:proof_univariate_exonential}

Let us introduce a family of functions $u_\lambda(x) = e^{\lambda(x - (\mu + t))}$. Then, for every $t \geqslant 0, \lambda \in \R$,
\begin{align*}
    \rho_{1}^{{\rm exp}}(\lambda, t) = \inf_{\alpha \in R, \beta \in \R} & \alpha + \beta \mu \ {\rm such \ that \ } \forall x \in X,  e^{\lambda(x - (\mu + t))} \leqslant \alpha + \beta x.
\end{align*}
The function $x \mapsto e^{\lambda(x - (\mu + t))} - (\alpha + \beta x)$ is convex on $[0, 1]$. Applying Proposition~\ref{prop:easy_exponential_cases}, we have:
\begin{align*}
    \rho_{1}^{{\rm exp}}(\lambda, t) = \inf_{\alpha \in R, \beta \in \R} \ \alpha + \beta \mu \ {\rm such \ that \ }  e^{-\lambda(\mu + t)} \leqslant \alpha \ {\rm and \ } e^{s(1 - (\mu + t))} \leqslant \alpha + \beta.
\end{align*}
The solution is given by $\alpha = e^{-\lambda(\mu + t)}$ and $\beta = e^{-\lambda(\mu + t)}(e^\lambda - 1)$, and thus we have $\rho_{1}^{{\rm exp}}(\lambda, t) = e^{-\lambda(\mu + t)}\left(1 + \mu(e^\lambda - 1) \right)$. By optimizing over $\lambda$, we have $\lambda_\star = \log\left(\frac{(\mu + t)(1 - \mu)}{\mu(1 - (\mu + t))}\right)$ and it holds for all $t \in [0, 1 - (\mu)]$ with $\nu = \mu + t$:
\begin{align*}
    \rho_{1, \star}^{\rm exp}(t) = \left(\frac{\mu}{\nu}\right)^\nu \left(\frac{1 - \mu}{1 - \nu}\right)^{1 - \nu}.
\end{align*}

\subsubsection{Alternative proof for Proposition~\ref{prop:univariate_exponential}}
\label{app:proof_exp_univariate_alternative}
Let $X$ be a random variable taking its value almost surely in $[0, 1]$, with mean $\E[X] = \mu$ and associated with the distribution $p \in \Px([0, 1])$. Then, for $KL(p, q) = \int_{X}\log(\frac{p(x)}{q(x)})dp(x)$ and $kl(\mu, \nu) = \log(\frac{\mu}{\nu})\mu$, the Kullback-Leibler divergence, it holds that
\begin{align*}
   \log\Prob(X \geqslant t) &\leqslant \inf_{s} -s t + \log(\E[e^{sX}]) {\rm \ by \ Markov's \ exponential \ inequality,}\\ 
   & \leqslant \inf_{s} -s t + \sup_{q} s\E_q[X] - KL(q, p), {\rm \ by \ Donsker-Varadhan's \ variational \ formula}, \\
   &= \inf_s -st + \sup_{\nu} s\nu - kl(\nu, \mu), \\
   &= - kl(t, \mu).
\end{align*}

\subsection{Variational approach with equal means: computing the extremal points for $\bar{X}_n$}
\label{app:extremal_points_equal_means}

Recall the optimization problem in the variational approach:
\begin{equation*}
\begin{aligned}
     \rho_n^{{\rm exp}} = \inf_{\alpha, \beta, t \geqslant 0} \  n(\alpha + \beta \mu - 1)
    {\rm \ such \ that \ } & \ -\sum_{i=1}^n \log(\alpha + \beta x_i^j) \leqslant 0, \ x^j \in {\rm extremal}(\bar{\X}_n), \\
    & \ \alpha + \beta x \geqslant 0, \ x \in [0, 1],
\end{aligned}
\end{equation*}
where $\bar{\X}_n = \left\{ (x_1, \ldots, x_n) \in [0, 1]^n,  x_1 + \cdots + x_n \geqslant n(\mu + t) \right\}$. Constraints $--\sum_{i=1}^n \log(\alpha + \beta x_i^j) \leqslant 0$ and $x_1 + \cdots + x_n \geqslant n(\mu + t)$ are symmetric in the coordinates, meaning that they have the same value when permuting $x_1^j, \ldots, x_n^j$. Therefore, the optimization problem requires to formulate only $O(n)$ extremal points, depending on the value for $n(\mu + t) \in [0, 1]$:
\begin{itemize}
    \item If $ n \leqslant n(\mu + t) > n-1$: ${\rm Extremal}(\X^n) = \{(1, \ldots, 1); (n(\mu + t) - (n-1), 1, \ldots, 1)\}$, 
    \item If $ n-2 < n(\mu + t) \leqslant n-1$: ${\rm Extremal}(\X^n) = \{(1, \ldots, 1); (0, 1, \ldots, 1); (0, n(\mu + t) - (n-2), 1, \ldots, 1)\}$, 
    \item If $ n-3 < n(\mu + t) \leqslant n-2$: ${\rm Extremal}(\X^n) = \{(1, \ldots, 1); (0, 1, \ldots, 1); (0, 0, 1, \ldots, 1); (0, 0, n(\mu + t) - (n-3), 1, \ldots, 1)\}$, 
    \item $\ldots$,
    \item If $ n(\mu + t) \leqslant 1$: ${\rm Extremal}(\X^n) = \{(1, \ldots, 1); (0, 1, \ldots, 1); \ldots; (0, \ldots, 0, 1); (0, \ldots, 0, n(\mu + t))\}$.
\end{itemize}
Thus, there are at most $\left\lfloor n(\mu + t)\right\rfloor + 1$ extremal points to consider.

\subsection{Closed-form solution to the variational reformulation for $n=2$}
\label{app:closed_form_reformulation_n2}

Recall the optimization problem under consideration for $n=2$, for two i.i.d. random variables taking their values in $[0, 1]$ with mean $\mu \in \R$,
\begin{equation}
\label{eq:variational_n2}
\begin{aligned}
    \log{ \rho_n^{{\rm exp}}}(t) = \inf_{\alpha, \beta, t \geqslant 0} \  2(\alpha + \beta \mu - 1)
    {\rm \ such \ that \ } & \ - \log(\alpha + \beta x_1^j) - \log(\alpha + \beta x_2^j) \leqslant 0, \ x^j \in {\rm extremal}(\bar{\X}_2),
\end{aligned}
\end{equation}
where $\bar{\X}_2 = \{(x_1, x_2) \in [0, 1]^2, x_1 + x_2 \geqslant 2(\mu + t)\}$. Extremal points of $\bar{\X}_2$ depends on the value of $\mu + t$:
\begin{itemize}
    \item If $ 1/2 \leqslant \mu + t \leqslant 1 $, ${\rm extremal}(\bar{\X}_2) = \{(1, 1), (1, 2(\mu + t) 1)\}$;
    \item If $0 \leqslant \mu + t \leqslant 1/2$,  ${\rm extremal}(\bar{\X}_2) = \{(1, 1),(0, 1), (0, 2(\mu + t)\}$.
\end{itemize}

We derive an optimal solution to problem~\eqref{eq:variational_n2} together with optimal values $(\alpha_\star, \beta_\star)$ in Proposition

\begin{proposition}
    Let $\mu \in [0, 1]$, and $t \in [0, 1 - \mu]$. Then, optimal solutions $(\alpha_\star, \beta_\star)$ to~\eqref{eq:variational_n2} verify: 
    \begin{align*}
        \alpha_\star &= \begin{cases} 
\sqrt{\frac{\mu}{\mu + t}} & \text{if } 0 \leqslant t \leqslant \frac{1}{2} - \mu, \\
\sqrt{\frac{1 - 2t - \mu}{1 - \mu}} - \frac{t(2(\mu + t) - 1)}{\sqrt{(1 - \mu)(1 - 2t - \mu)}(1 - (\mu + t))} & \text{if } \frac{1}{2} - \mu \leqslant t \leqslant \frac{(1 - \mu)^2}{2 - \mu}, \\
0 & \text{if } t \geqslant \frac{(1 - \mu)^2}{2 - \mu},
\end{cases} \\
     \beta_\star &= \begin{cases} 
\frac{t}{(\mu + t)\sqrt{\mu(\mu + 2t)}} & \text{if }  0 \leqslant t \leqslant \frac{1}{2} - \mu, \\
\frac{t}{(1 - (\mu + t)) \sqrt{(1 - \mu)(1 - \mu - 2t)}} & \text{if } \frac{1}{2} - \mu \leqslant t \leqslant \frac{(1 - \mu)^2}{2 - \mu}, \\
\frac{\mu}{\sqrt{2(\mu + t) - 1}} & \text{if } t \geqslant \frac{(1 - \mu)^2}{2 - \mu}.
\end{cases}
    \end{align*}
\end{proposition}
\textbf{Proof. } The proof is divided into two parts, depending on the value for $\mu + t$.

1. Let us first assume that $\mu + t \in ]1/2, 1]$. Then, the optimization problem~\eqref{eq:variational_n2} takes the form:
\begin{align*}
    \inf_{\alpha, \beta} \ 2(\alpha + \beta \mu - 1), \ {\rm such \ that \ }& \alpha \geqslant 0, \ &(\times \lambda_1) \\
    & -2\log(\alpha + \beta) \leqslant 0,  \ &(\times \lambda_2)\\
    & - \log(\alpha + \beta ) - \log(\alpha + (2(\mu + t) - 1) \beta) \leqslant 0 \ &(\times \lambda_3), \\
\end{align*}
where $\lambda_1, \lambda_2, \lambda_3 \geqslant 0$ are dual variables. Reformulating the second constraint into ``$\alpha + \beta \geqslant 1$'', it still holds that this problem is a convex. We compute the KKT conditions:
\begin{align*}
    2 &= \lambda_1 + \lambda_2 + \lambda_3\left( \frac{1}{\alpha + \beta} + \frac{1}{\alpha + (2(\mu + t) - 1) \beta}\right), \\
    2 \mu &= \lambda_2 + \lambda_3 \left(\frac{1}{\alpha + \beta} + \frac{ (2(\mu + t) - 1)}{\alpha + (2(\mu + t) - 1) \beta} \right), \\
    0 &= \lambda_1 \alpha , \\
    0 &=\lambda_2 (1 - (\alpha + \beta)), \\
    0 &= \lambda_3 (- \log(\alpha + \beta ) - \log(\alpha + (2(\mu + t) - 1) \beta) ), \\
   0 & \leqslant \lambda_1, \lambda_2, \lambda_3.
\end{align*}
We proceed by a distinction of cases:
\begin{itemize}
    \item Assume first that $\lambda_1 \neq 0$, then $\alpha = 0$. If $\lambda_2 \neq 0$, then $\beta = 1$ and $\lambda_3 =  0 $, $\lambda_2 = 2 = 2\mu$, which is false. We conclude that $\lambda_2 = 0$, and that $\lambda_3 \neq 0$. This fixes the value for $\beta$, by solving $\log(\beta) + \log(\beta (2 (\mu + t) - 1)) = 0$ and injecting this relationship into the two first KKT conditions above.
    \item Assume now that $\alpha \neq 0$. Then, $\lambda_1 = 0$. At least $\lambda_2$ or $\lambda_3$ must be nonzero, and both cannot be nonzero if $t + \mu \neq 1/2$. Assuming $ \log(\alpha + \beta ) + \log(\alpha + (2(\mu + t) - 1) \beta) \neq 0$ entails $\lambda_3 = 0$ and $\lambda_2 \neq 0$, that is $\alpha + \beta = 1$. Then, $\log(\alpha + (2(\mu + t) - 1) \beta)< 0$ which is false. We conclude that $ \log(\alpha + \beta ) + \log(\alpha + (2(\mu + t) - 1) \beta) = 0$, and thus, that $\lambda_3 \neq 0$ and $\lambda_2 = 0$. This leads to the solutions : 
    \begin{align*}
        \alpha_\star &= \sqrt{\frac{1 - 2t - \mu}{1 - \mu}} - \frac{t(2(\mu + t) - 1)}{\sqrt{(1 - \mu)(1 - 2t - \mu)}(1 - (\mu + t))}, \\
        \beta_\star&= \frac{t}{(1 - (\mu + t)) \sqrt{(1 - \mu)(1 - \mu - 2t)}},
    \end{align*}
    for which $\alpha_\star > 0$ for $t \leqslant \frac{(1 - \mu)^2}{2 - \mu}$.
\end{itemize}

2. We now consider the case where $\mu + t \in [0, 1/2]$. The optimization problem under consideration takes the form:
\begin{align*}
    \inf_{\alpha, \beta}  \ 2(\alpha + \beta \mu - 1), \ {\rm such \ that \ }& \alpha \geqslant 0, \ &(\times \lambda_1) \\
    & \alpha + \beta \geqslant 1,  \ &(\times \lambda_2)\\
    & - \log(\alpha) - \log(\alpha + \beta) \leqslant 0 \ &(\times \lambda_3), \\
    & - \log(\alpha) - \log(\alpha + 2(\mu + t) \beta) \leqslant 0 \ &(\times \lambda_4).
\end{align*}
We compute the KKT conditions, leading to
\begin{align*}
    2 &= \lambda_1 + \lambda_2 + \lambda_3\left(\frac{1}{\alpha} + \frac{1}{\alpha + \beta}\right) + \lambda_4 \left(\frac{1}{\alpha} + \frac{1}{\alpha + \beta 2(\mu + t)} \right), \\
    2 \mu &= \lambda_2 + \lambda_3 \frac{1}{\alpha + \beta} + \lambda_4 \frac{2(\mu + t)}{\alpha + \beta 2(\mu + t)}, \\
   0& =  \lambda_1 \alpha , \\
    0& =\lambda_2 (1 - (\alpha + \beta)), \\
    0& =\lambda_3 (- \log(\alpha) - \log(\alpha + \beta) ), \\
   0& = \lambda_4 (- \log(\alpha) - \log(\alpha + \beta2(\mu + t))), \\
   0 &\leqslant \lambda_1, \lambda_2, \lambda_3, \lambda_4 .
\end{align*}
First, note that $\alpha > 0$, and thus, $\lambda_1 = 0$. In addition, $\lambda_2 = 0$, otherwise, it implies that $\lambda_3 = \lambda_4 = 0$ and $\mu = 1$, which is impossible. In addition, if $\lambda_3 \neq 0$, then $\alpha(\alpha + \beta) = 1$. Yet, since $\alpha > 0$, $\beta \geqslant 0$. The function $t \mapsto \alpha(\alpha + t \beta)$ is nondecreasing, and thus $\alpha (\alpha + 2(\mu + t)\beta) > 1$, which is impossible. We conclude that $\lambda_3 = 0$. 
By construction, $\lambda_4 \neq 0$, which implies that $\alpha (\alpha + 2 (\mu + t)\beta) = 1$, and thus,
\begin{align*}
        2 &= \lambda_4(\frac{1}{\alpha} + \alpha), \\
        2 \mu &= \lambda_4 \alpha 2 (\mu + t),
\end{align*}
from which we conclude the solutions $\alpha = \sqrt{\frac{\mu}{\mu + t}}$ and $\beta = \frac{t}{(\mu + t)\sqrt{\mu(\mu + 2t)}}$.
\hfill \Halmos

\subsection{Proof for Lemma~\ref{lem:nbconstraints_two_blocks}.}
\label{app:proof_lemma_constraints_two_blocks}
Problem~\eqref{eq:var_diff_mean} simplifies into 
\begin{align*}
    \inf_{\alpha \in \R^2, \beta \in \R^2} & m(\alpha_1 + \beta_1 \mu_1) + (n - m)(\alpha_2 + \beta_2 \mu_n), \\
    {\rm \ such \ that \ } & \forall x \in {\rm extremal(\bar{\X}_n}), \ -\sum_{i=1}^m\log(\alpha_1 + \beta_1 x_i) -\sum_{i=m+1}^n \log(\alpha_2 + \beta_2 x_i) \leqslant 0,
\end{align*} 
where $\bar{\X}_n = \left\{(x_1, \ldots, x_n) \in [0, 1]^n, x_1 + \ldots+ x_n\geqslant nt + m\mu_1 + (n - m)\mu_n \right\}$. We define $\bar{\mu}_n = \frac{1}{n}\sum_{i=1}^n \mu_i$ and denote $q = n(\bar{\mu}_n + t) - \lfloor n(\bar{\mu}_n + t) \rfloor$. We denote the points with the notation $(x_1, \ldots x_m | x_{m+1}, \ldots, x_n)$. Let $n(\bar{\mu}_n + t) \in [n-k+1, n-k[$. The following assertions are true:
\begin{itemize}
    \item $(1, \ldots, 1) \in {\rm extremal}(\bar{\X}_n)$
    \item If $(1, \ldots, 1|0, 1, \ldots, 1) \in {\rm extremal}(\bar{\X}_n)$, all its permutations are in ${\rm extremal}(\bar{\X}_n)$. By symmetry, it is sufficient to notice that the point $(1, \ldots, 1, 0| 1, \ldots, 1)$ is in ${\rm extremal}(\bar{\X}_n)$.
    \item If $(1, \ldots, 1|0, 0, 1, \ldots, 1) \in {\rm extremal}(\bar{\X}_n)$, all its permutations are in ${\rm extremal}(\bar{\X}_n)$. By symmetry, it is sufficient to notice that the point $(1, \ldots, 1, 0| 0, 1, \ldots, 1)$ and $(1, \ldots, 1, 0, 0| 1, \ldots, 1)$ are in ${\rm extremal}(\bar{\X}_n)$.
    \item If $(1, \ldots, 1|0, \ldots 0, 1, \ldots, 1) \in {\rm extremal}(\bar{\X}_n)$, all its permutations are in ${\rm extremal}(\bar{\X}_n)$. By symmetry, it is sufficient to notice that the point $(1, \ldots, 1, 0| 0, \ldots, 0, 1, \ldots, 1), \ldots, (1, \ldots, 1, 0, \ldots, 0| 1, \ldots, 1)$ are in ${\rm extremal}(\bar{\X}_n)$ (as long as $k \geqslant m$). That is about $O(|k-m|)$ points.
    \item If $(1, \ldots, 1| 0, \ldots, 0, q, 1, \ldots, 1) \in {\rm extremal}(\bar{\X}_n)$, all its permutations are in ${\rm extremal}(\bar{\X}_n)$. By symmetry, it is sufficient to notice that points $(1, \ldots, 1, q|0, \ldots, 0, 1, \ldots, 1)$, $(1, \ldots, 1, q, 0 | 0,\ldots, 0, 1, \ldots, 1)$, $\ldots$, $(1, \ldots, 1, q, 0, \ldots, 0 | 1, \ldots, 1)$ are in ${\rm extremal}(\bar{\X}_n)$.  That is, about $O(2|k-m|)$ points.
\end{itemize} 
If they are $k$ zeros elements, there are $O(\sum_{i=1}^{k - m} i )= O(k^2) \leqslant O(n^2)$ points.

\section{Connecting the second-order to the first-order moment: proof for Lemma~\ref{lem:var_to_mean_relationship}}
\label{proof_lem:var_to_mean_rel}
Let $X$ be a random variable almost surely in $[0, 1]$ with mean $\mu^{(1)}$. First, its second-order moment exists. In addition, we have, 
\begin{align*}
    \mu^{(2)} &\leqslant \sup_{p \in \mathcal{P}([0, 1])} \int_{0}^1 x^2 dp(x), {\rm \ such \ that \ } \int_{0}^1 x dp(x) = \mu^{(1)}, \\
    & = \inf_{\alpha, \beta} \alpha + \beta \mu^{(1)} {\rm \ such \ that \ } \forall x \in [0, 1], x^2 \leqslant \alpha + \beta x, \\
    &= \inf_{\alpha, \beta} \alpha + \beta \mu^{(1)} {\rm \ such \ that \ } \alpha \geqslant 0, 1 \leqslant \alpha + \beta {\rm \ \ (by \ convexity), \ } \\
    &= \mu^{(1)},
\end{align*}
with $\alpha = 0, \beta =1$. In addition, by Cauchy-Schwarz, $\int_0^1 x dp(x) \leqslant \sqrt{\left( \int_0^1 x^2 dp(x)\right)}$ that is $\mu^{(1)} \leqslant \sqrt{\mu^{(2)}}$.

\end{appendices}

\ACKNOWLEDGMENT{This work was funded by MTE and the Agence Nationale de la Recherche as part of the “Investissements d’avenir” program, reference ANR-19-P3IA-0001 (PRAIRIE 3IA Institute). We also acknowledge support from the European Research Council (grant SEQUOIA 724063).}	

\bibliographystyle{ormsv080} 
\bibliography{references.bib} 

\end{document}